\documentclass[a4paper,11pt]{article}
\usepackage{amsfonts}
\usepackage{amssymb}
\usepackage{amsmath}
\usepackage{amsthm}

\usepackage{multirow}
\usepackage{diagbox}
\usepackage{amsopn}
\usepackage{gensymb} 
\usepackage{fancyhdr}
\usepackage{graphicx}
\usepackage{mathrsfs}
\usepackage{latexsym}
\usepackage{graphicx}
\usepackage{subfigure}

\usepackage{subfig}
\usepackage{caption}

\usepackage{color}
 \usepackage{calc}
  \usepackage{tikz}
\usepackage{etoolbox} 
\usepackage{listofitems} 

\tikzset{>=latex} 
\colorlet{myred}{red!80!black}
\colorlet{myblue}{blue!80!black}
\colorlet{mygreen}{green!60!black}
\colorlet{mydarkred}{myred!40!black}
\colorlet{mydarkblue}{myblue!40!black}
\colorlet{mydarkgreen}{mygreen!40!black}
\tikzstyle{node}=[very thick,circle,draw=myblue,minimum size=22,inner sep=0.5,outer sep=0.6]
\tikzstyle{connect}=[->,thick,mydarkblue,shorten >=1]
\tikzset{ 
  node 1/.style={node,mydarkgreen,draw=mygreen,fill=mygreen!25},
  node 2/.style={node,mydarkblue,draw=myblue,fill=myblue!20},
  node 3/.style={node,mydarkred,draw=myred,fill=myred!20},
}

\usepackage{setspace}
\usepackage{textcomp} 
\usepackage{neuralnetwork}

 \usepackage[noend]{algpseudocode}

\usepackage{algorithmicx,algorithm}
 
 \usepackage{appendix}
 
 \usepackage{graphicx,multicol}
 
\usetikzlibrary{decorations.pathreplacing}

\numberwithin{equation}{section} \allowdisplaybreaks
\theoremstyle{plain}

\theoremstyle{remark}
\newtheorem{remark}{Remark}
\numberwithin{equation}{section}

\begin{document}


\title{\bf{Bayesian approach for limited-aperture inverse acoustic scattering with total variation prior}}

\author{
\hspace{-2.5cm}{\small    Xiao-Mei Yang$^1$, Zhi-Liang Deng$^{2}$\thanks{
Corresponding author: dengzhl@uestc.edu.cn;
Supported by Central Government Funds of Guiding Local Scientific and Technological Development for Sichuan Province No. 2021ZYD0007, NSFC No. 11601067.}  and Ailin Qian$^{3}$ } \\
\hspace{-1.5cm} {\scriptsize
$1.$ School of Mathematics,
Southwest Jiaotong University,
Chengdu 610031, China}
 \\
\hspace{-1.5cm}{\scriptsize
$2.$ School of Mathematical Sciences,  University of Electronic Science and Technology of China,
Chengdu 610054, China
}
 \\
\hspace{-1.5cm}{\scriptsize
$3.$ School of Mathematical Sciences, Hubei University of Science and Technology, Xianning 437100, China
}
}
\date{}
\maketitle

\begin{abstract}
\noindent 
 In this work, we apply the Bayesian approach for the acoustic scattering problem to reconstruct the shape of a sound-soft obstacle using the limited-aperture far-field measure data. 
A novel total variation prior is  assigned to the shape parameterization form.
This prior is imposed on the Fourier coefficients of the parameterized form of the obstacle.
Extensive numerical tests are provided  to illustrate the numerical performance.


\noindent \textbf{Key words:}    Inverse acoustic scattering; Uncertainty quantification; $l_1$ priors

\noindent \textbf{MSC 2010}: 65R32, 65R20
\end{abstract}

\section{Introduction}

 This investigation is concerned with the 2-dimensional acoustic obstacle reconstruction problem. 
 The goal is to detect and identify the unknown object using acoustic wave measured data.
 Such problem  is of practical importance  and   has attracted a lot of interest.
 Over the past years, many satisfactory reconstructions have been provided \cite{Cakoni, Colton1, Kirsch}.
 For more references, one can refer to \cite{Colton1} for linear sampling method,  \cite{Kirsch, Kirsch1} for the factorization method,  \cite{Ito, JLi} for the direct sampling method and  \cite{Bui, Li} for Bayesian approach.  
 As we can see, most of these methods focus on the reconstruction problem using the full-aperture data. The limited-aperture case usually arises for example from practical difficulties in surrounding the whole region of interest by sensors.   Some methods to full-aperture data and some new ideas have been developed to process the limited-aperture case, e.g., \cite{ Ito,   Li,  Liuj, Liux, Zinn, Potthast1}.  
Since lack of full aperture measurements, the ill-posedness and nonlinearity of the inverse problem become more severe, which makes this problem more difficult to solve. The reconstruction is not as good as the full aperture case in general \cite{Audibert, Guo}.
 
 

The present work aims at the application of the Bayesian statistical approach in the obstacle reconstruction problem, which is a further development of \cite{Li}.
Our contribution in this paper is to propose a novel prior form for the unknown shape parameterization.
This prior is a modified version of total variation prior. We still call it as the total variation prior (TV). Specifically, we suppose that the shape parameterization representation is expanded in Fourier series. The Fourier coefficients are assumed in the total variation prior.
In \cite{Bui, Li}, the prior is imposed on the parameterization representation using Gaussian random field prior. Subsequently, this random field prior can be considered as a mutually independent Gaussian prior on Fourier coefficients by virtue of Karhunen-Lo\`eve expansion.
 In \cite{Kaipio},  we can see that the TV Markov random field prior admits advantages in designing structure priors. It is generally defined by the way of neighbor system. The TV prior in this paper is generated from the Gaussian random field prior \cite{Bui, Li} by the idea in \cite{Wang}. 
 The induced parameterization representation is not as smooth as the representation from the Gaussian random field. 



In the remainder, this paper is organized as follows:  In Section 2, we introduce the inverse acoustic obstacle scattering problems and apply the Bayesian approach to solve it. In Section 3, a new TV prior is introduced. In Section 4, we  give some numerical examples to illustrate the numerical effectiveness. Finally, the conclusions is stated. 


\section{Inverse Scattering Problems and Bayesian approach}
We begin with the formulations of the acoustic scattering problem. 
Let $\Omega \subset \mathbb{R}^{2}$ be a bounded, simply connected domain with $C^{2}$ boundary $\partial \Omega$. Define $\mathbb{S}=\left\{x \in \mathbb{R}^{2},|x|=1\right\}$. 
Consider a  time-harmonic plane wave propagating in the exterior of $D$. 
The velocity potential is of this form
\begin{align}
U(x, t)=\operatorname{Re}\left\{u^{\rm i}(x){\mathrm e}^{-i\omega t}\right\} \,\,\text{with}\,\, u^{\rm i}(x):=\mathrm{e}^{i \kappa x \cdot d}, \,\, x \in \mathbb{R}^{2}, \,\, d \in \mathbb{S},
\end{align}
where $d$ is the propagating direction,  $\omega>0$ is frequency, $\kappa=\omega/c>0$ is the wavenumber and $c$ is the speed of sound. Hereafter,  the term $\exp(-i\omega t)$ is factored out and only the space dependent part of all waves is considered.
Let $D$ be an impenetrable sound-soft obstacle. 
Clearly, the obstacle will "scatter" the wave $u^{\rm i}$ so that we may write the total field as $u=u^{\rm i}+u^{\rm s}$, where $u^{\rm s}$ is the scattered wave. 
The classical physical problem is to explore the scattering phenomenon, i.e., find $u^{\rm s}$. As we know, the scattered field $u^{\rm s}$ is governed by an exterior boundary value problem for Helmholtz equation
\begin{subequations}\label{eqn0}
\begin{align}
&\triangle u+\kappa^{2} u =0, &\text { in } &\mathbb{R}^{2} \backslash \bar{\Omega},\label{eqn1}\\
&u =0, &\text { on }& \partial \Omega,\label{eqn2} \\
&\lim _{r \rightarrow \infty} \sqrt{r}\left(\frac{\partial u^{\rm s}}{\partial r}-i \kappa u^{\rm s}\right) =0, &&\label{eqn3}
\end{align}
\end{subequations}
where  equation \eqref{eqn2} is the sound-soft boundary condition and \eqref{eqn3} is the Sommerfeld radiation condition.

The well-posedness of the direct scattering problem \eqref{eqn1}-\eqref{eqn3} has been
established, see, e.g., \cite{Colton}. We refer to the standard monograph \cite{Colton} for a research statement on the significant progress both in the mathematical theories
and the numerical approaches.
It is well-known that the unique solution $u^{\rm s}$ of \eqref{eqn0}  admits an asymptotic expansion \cite{Colton}
\begin{align}
u^{\rm s}(x, d)=\frac{\mathrm{e}^{i \frac{\pi}{4}}}{\sqrt{8 \kappa \pi}} \frac{\mathrm{e}^{i \kappa r}}{\sqrt{r}}\left\{u^{\infty}(\hat{x}, d)+\mathcal{O}\left(\frac{1}{r}\right)\right\} \quad \text { as } r:=|x| \rightarrow \infty
\end{align}
uniformly in all directions $\hat{x}=x /|x|$. The function $u^{\infty}(\hat{x}, d)$ is called the far-field pattern. The  inverse scattering problem in this paper is to find the shape of obstacle $\Omega$ from some observed far-field data.
Let $\gamma^{\rm o}\subseteq(\subsetneq) \mathbb{S}$ be the observation aperture and $\gamma^{\rm i}\subsetneq \mathbb{S}$ be the incident direction.  
The direct scattering problem can be formulated as
\begin{align}\label{par1}
u^{\infty}(\hat{x}, d)=\mathcal{F}(\Omega), \quad(\hat{x}, d) \in \gamma^{\rm o} \times \gamma^{\rm i}
\end{align}
where $\mathcal{F}$ is the shape-to-measurement operator. If $\gamma^{\rm o}\subsetneq\mathbb{S}$, i.e., $\gamma^{\rm o}$ is a proper subset of $\mathbb{S}$,  this case  is referred as the limited-aperture problem. 
If $\gamma^{\rm o}=\mathbb{S}$, we are concerned with the full-aperture case.

We consider a starlike obstacle $\Omega$ centering at $(x_{\rm c}, y_{\rm c})$ with boundary being parametrized as
\begin{align}
\partial \Omega:=(x_{\rm c}, y_{\rm c})+r(\theta)(\cos \theta, \sin \theta)=(x_{\rm c}, y_{\rm c})+\exp (q(\theta))(\cos \theta, \sin \theta),
\end{align}
where $q(\theta)=\log r(\theta),\,\, 0<r(\theta)<r_{\max }$, $\theta \in[0,2 \pi)$. With the parameterization, we take the noise in measurements into account and rewrite  \eqref{par1} as a statistical inference model
\begin{align}\label{mod1}
y=\mathcal{F}(q)+\eta,
\end{align}
where $q \in X$ and $y=u_{\infty}(\hat{x}, d) \in Y$ for some suitable Banach spaces $X$ and $Y$, respectively. In particular, $y$ is the noisy observations of $u^{\infty}(\hat{x}, d)$ and $\eta(\hat{x}, d)$ is the noise. In this paper, we assume that the observation noise is normal with mean zero and independent of $q$, i.e., $\eta(\hat{x}, d) \sim \mathcal{N}\left(0, \sigma^{2}(\hat{x})\right)$.

Apart from mere estimation purposes, a major goal of statistical inference is to find the posterior distribution $p(q|y)$.  Bayes Theorem provides a principled way for calculating the posterior conditional probability.
According to Bayes' formula, it has the following representation
\begin{align}\label{in2.4}
p(q| y)=\frac{p(y|q)p(q)}{p(y)}\propto p(y|q)p(q),
\end{align}
where $p(y|q)$ is the likelihood, $p(q)$ is the prior and $p(y)$ is the normalized constant. The likelihood model is established by the forward model and noise model. And it admits by \eqref{mod1}
\begin{align}\label{in2.5}
p(y|q)\propto\exp\left(-\frac{\|\mathcal{F}(q)-y\|^2}{2\sigma^2}\right):=\exp\left(-\Phi(q, y)\right).
\end{align}
The prior model provides the prior knowledge for the unknown parameters before the data is collected. We will discuss it in the subsequent section.
One also can refer to \cite{Bardsley, Calvetti, Stuart} for more details. 

The Bayesian approach answers richer questions than the traditional regularization-based methodologies. As we can see, the solution it produces is not limited to individual values but consists of probability distributions.  
Therefore, it provides an effective tool to analyze the relative probabilities of different approximate solutions or the probability of a solution to lie in a subset of the solution space for decision maker. 


\section{Total variation prior}

A key sticking point of Bayesian analysis is the prior elicitation, which usually plays a defining role in Bayesian statistical inference.  There is a vast literature on the related discussion, e.g.,  \cite{Ali, Kaipio, Wang}.
A rich class of prior distributions can be derived from the theory of Markov random fields (MRF), which is usually defined in the discrete form by the neighborhood system \cite{Bardsley, Kaipio, Wang}. 
In this paper, the prior for $q$ is taken as the Gauss MRF $N(0, \Gamma)$ with covariance operator $\Gamma=Q^{-1}$. Here $Q$ is the precision operator defined by the  fractional diffusion operator $-(\frac{d^2}{d\theta^2})^{s}$ with periodic boundary condition. 

In real numerical implementation, we need the discrete form of a random variate. There are two ways to transform the infinite dimensional random variate into its finite pattern. 
In the first way, we take a mesh partition $\tilde{\Omega}$. Denote the discretized version of precision operator $Q$ by $\operatorname{Q}$. The discretize can be done by the finite difference, finite element etc.
 The joint density function for $\boldsymbol{q}:=[q(\theta_1), q(\theta_2), \cdots, q(\theta_m)]^\top$ is Gaussian and admits the form
\begin{align}\label{dis1}
p(\boldsymbol{q})=(2\pi)^{-m/2}|\operatorname{Q}|^{1/2}\exp(-\frac{1}{2}\boldsymbol{q}^\top \operatorname{Q} \boldsymbol{q}).
\end{align}
In this way, the finite dimensional samples of $q$ can be generated by the probability distribution \eqref{dis1}. Denote eigenvalues and eigenvectors of precision matrix $\operatorname{Q}$ by $\lambda_{\operatorname{Q}}$ and $q_{\operatorname{Q}}$.
In the second way, we use the Karhunen-Lo\`eve expansion to represent the prior random variate. The precision operator $Q$ admits eigenvalues $\lambda_k$ and eigenfunctions $\psi_k$, $k=1, 2, \cdots$, i.e.,
\begin{align}
Q \psi_k=\lambda_k \psi_k, \,\, k=1, 2, \cdots.
\end{align}
According to Karhunen-Lo\`eve expansion of $q$, we have
\begin{align}\label{kal1}
q(x)=\sum_{k=1}^\infty \frac{B_k}{\lambda_k}\psi_k(x)\approx\sum_{k=1}^m \frac{B_k}{\lambda_k}\psi_k(x):=\tilde{q}(x),
\end{align}
where $B_k$s are independent identity distribution normal Gaussian. The truncated form $\tilde{q}(x)$ is used in the numerical implement process.

In fact, it can be verified that $\lambda_{\operatorname{Q}}$ behaves like $\lambda_k$ and $\theta_{\operatorname{Q}}$s  are vectors of values of $\psi_k$ at discrete points. Based on this reason, in our numerical test, we do not distinguish them any more. 
And in the subsequent parts, we only use the KL expansion to represent the prior. Therefore, the prior for the unknown $q$ is determined by its expansion coefficients $B_k$s. 
\begin{remark}
In the KL expansion prior condition, we actually assume the  function $q(\theta)$ admits Fourier series representation
\begin{align}
q(\theta)=\sum_{k=1}^\infty B_k \psi_k(\theta).
\end{align}
And the prior for $B_k$ is given by $N(0, \frac{1}{\lambda_k^2}I)$ from \eqref{kal1}. 
\end{remark}

For a Gaussian random vector with mutually independent components, we may transform it to an $l_1$-type prior according to the theories in \cite{Wang}.  
For this purpose, an invertible mapping function $g: \mathbb{R} \rightarrow \mathbb{R}$ is introduced
\begin{align}
g(B_k) \equiv \mathcal{L}^{-1}(\mathcal{G}(B_k))=-\frac{1}{\lambda} \operatorname{sign}(B_k) \log \left(1-|2\mathcal{G}(B_k)-1|\right),
\end{align}
where $\mathcal{L}$ is the cumulative distribution function (cdf) of the Laplace distribution and $\mathcal{G}$ is the cdf of the standard Gaussian distribution. This function $g$ relates a Gaussian reference random variable $B_k \in \mathbb{R}$ to the Laplace-distributed parameter $A_k \in \mathbb{R}$, such that $A_k=g(B_k)$. Then we apply the transformation to each components of $B$ 
\begin{align}\label{prior1}
A(B):=[g(B_1), g(B_2), \cdots, g(B_m)]^\top.
\end{align}
 A prior transformation for the $l_1$-type prior is $Z_q:=D^{-1}A(B)$, i.e.,
 \begin{align}\label{prior2}
 p(Z_q)\propto \exp(-\lambda |DZ_q|)=\exp(-\lambda \sum_{i=1}^m |DZ_q|_i),
 \end{align}
 where $D$ is an invertible matrix. For total variation prior, the invertible matrix $D$ is taken as \cite{Wang}
 \begin{equation}\label{prior3}
D_0=\left[\begin{array}{cccc}
1 & & & 1 \\
-1 & 1 & & \\
& \ddots & \ddots & \\
& & -1 & 1
\end{array}\right]_{m \times m}.
\end{equation}
 We use a modified version of $D_0$ as $D=I+\alpha D_0$ with a small constant $\alpha_0>\alpha>0$.  By this, we actually take $|DZ_q|$ in \eqref{prior2} as
 \begin{align}\label{nons1}
 |DZ_q|=|(1+\alpha)z_1+\alpha z_m|+\sum_{i=2}^m |z_i+\alpha(z_i-z_{i-1})|,
 \end{align}
 where $z_i$ is the $i$-th component of $Z_q$. 
 Since $D$ is nonsingular matrix, \eqref{nons1} defines a norm of vector $Z_q$.  Therefore, it follows that the norm in \eqref{nons1} is equivalent to the normal norm. 
\begin{remark}
The total variation strengthens the relatedness of the Fourier coefficients. In \eqref{kal1}, the Fourier coefficients are required to be independent. However, the convergence of Fourier series actually implies that they are not completely uncorrelated. Though we here can not use the total variation prior to reflect the convergence, the correlation of Fourier coefficients is added. 
\end{remark}


\section{Numerical examples}
We use the pCN Markov chain Monte Carlo algorithm \cite{Stuart} listed in Algorithm \ref{alg1} to implement the sampling process. 
\begin{algorithm}
\caption{pCN MCMC algorithm}
\begin{algorithmic}[1]

    \State Initialize: Pick $B^0$. Promote $B^0$ to $Z^0_q$ according to \eqref{prior1}-\eqref{prior3} and compute $\tilde{q}(Z^0_q)$ by \eqref{kal1}. Set $j=0$.
    \State While $j<N$, continue: 
    \begin{itemize}
\item Draw a proposal $\hat{B}$ according to
\begin{align*}
\hat{B}=\sqrt{1-\beta^2}B^j+\beta B_{\rm pr},
\end{align*}
where $B_{\rm pr}\sim N(0, I)$ and $\beta\in (0, 1)$ is a constant.
Transform $\hat{B}$ to $\hat{Z}_q$ according to \eqref{prior1}-\eqref{prior3}. Generate $\tilde{q}(\hat{Z}_q)$ using \eqref{kal1} with coefficients $\hat{Z}_q$.
\item Compute the acceptance probability:
\begin{align*}
\alpha=\min\left\{1, \frac{\exp\left(-\Phi(\tilde{q}({\hat{Z}_q}), g)\right)}{\exp\left(-\Phi(\tilde{q}(Z_q^j), g)\right)}\right\}.
\end{align*}
\item Flip an $\alpha$-coin: Draw $\xi\sim \operatorname{Uniform}([0, 1])$:
\begin{itemize}
\item If $\xi\leq \alpha$, accept the proposal, setting $B^{j+1}=\hat{B}$,
\item else, reject the move and stay put $B^{j+1}=B^j$.
\end{itemize}
\item Increase $j\leftarrow j+1$. 
\end{itemize}

\end{algorithmic}
\label{alg1}
\end{algorithm}
We take several benchmark examples in acoustic obstacle scattering to verify the numerical effectiveness ($\theta \in(0,2 \pi]$), see Table \ref{table:nonlin}.
\begin{table}[ht]
\caption{The boundary parameterizations of the object obstacles} 
\centering 
\small
\begin{tabular}{l l  l} 
\hline 
Kite &   & $\left(x_{1}, \, x_{2}\right)=(\cos \theta+0.65 \cos 2 \theta-0.65,\,1.5 \sin \theta)$ \\ [1ex]\\ 
Roundrect &   & $r(\theta)=\left(\cos^4\theta+(2/3\sin\theta)^4\right)^{-1/4}$ \\ [1ex]\\
Acorn &   & $r(\theta)=\frac{3}{5}\sqrt{\frac{17}{4}+2\cos3\theta}$ \\ [1ex]\\
Pear &  & $r(\theta)=\frac{5+\sin 3 \theta}{6}$ \\ [1ex]\\
Bean &  & $r(\theta)=0.4\sqrt{4\cos^2\theta+\sin^2\theta}$\\ [1ex] \\
Threelobes & &$ r(\theta)=0.5+0.25\exp(-\sin3\theta)-0.1\sin\theta$\\ [1ex]\\
Star && $r(\theta)=1+0.3\sin5\theta$\\ [1ex]\\
Cloverleaf && $r(\theta)=1+0.3\cos4\theta$\\[1ex]\\
Peanut && $r(\theta)=0.4\sqrt{4\cos^2\theta+\sin^2\theta}$\\[1ex] \\ 
Drop && $\left(x_1,\, x_2\right)=(-1+2\sin(\theta/2),\, -\sin\theta)$\\
\hline 
\end{tabular}
\label{table:nonlin} 
\end{table}

The true data is generated by solving the forward problem by boundary integral equation method \cite{Colton}. 
Let $\phi$  be the observation angle such that $\hat{x}:= (\cos \phi, \sin \phi )$. We use the  observation/measurement apertures as follows:
\begin{equation*}
\begin{aligned}
&\gamma_{1}^{\rm o}=\{(\cos \phi, \sin \phi) \mid \phi \in[0,2 \pi]\}, \\
&\gamma_{2}^{\rm o}=\{(\cos \phi, \sin \phi) \mid \phi \in[0, \pi]\},\\
&\gamma_{3}^{\rm o}=\{(\cos \phi, \sin \phi) \mid \phi \in[0, \pi / 2]\}.
\end{aligned}
\end{equation*}
The incident apertures are
$$
\begin{aligned}
\gamma_{1}^{\rm i} &=\{(1,0)\}, \\
\gamma_{2}^{\rm i} &=\left\{(\cos \theta, \sin \theta) \mid \theta=\{\pi / 2, 3 \pi / 2\}\right\}.
\end{aligned}
$$
The relative noise is added to the true data by 
\begin{align}
y=u^\infty(\hat{x}, d)+(\eta_1+\eta_2 i)\|u^\infty(\hat{x}, d)\|_\infty.
\end{align}

In these numerical tests, we fix some parameters in the prior, $s=2.2$, $\lambda=0.2$, $\alpha=0.1$. The Fourier expansion of $q$ is truncated to the first $27$ terms. The wave number is set to $\kappa=1$. The true center position $(x_{\rm c}, y_{\rm c})$ is fixed at $(0, 0)$. And we take the samples after $1000$ iteration and use the mean estimation as the approximation of the unknown parameters. 

\begin{figure}
\captionsetup[subfigure]{labelformat=empty}
    \centering
    \subfigure[Kite]{\includegraphics[width=4cm]{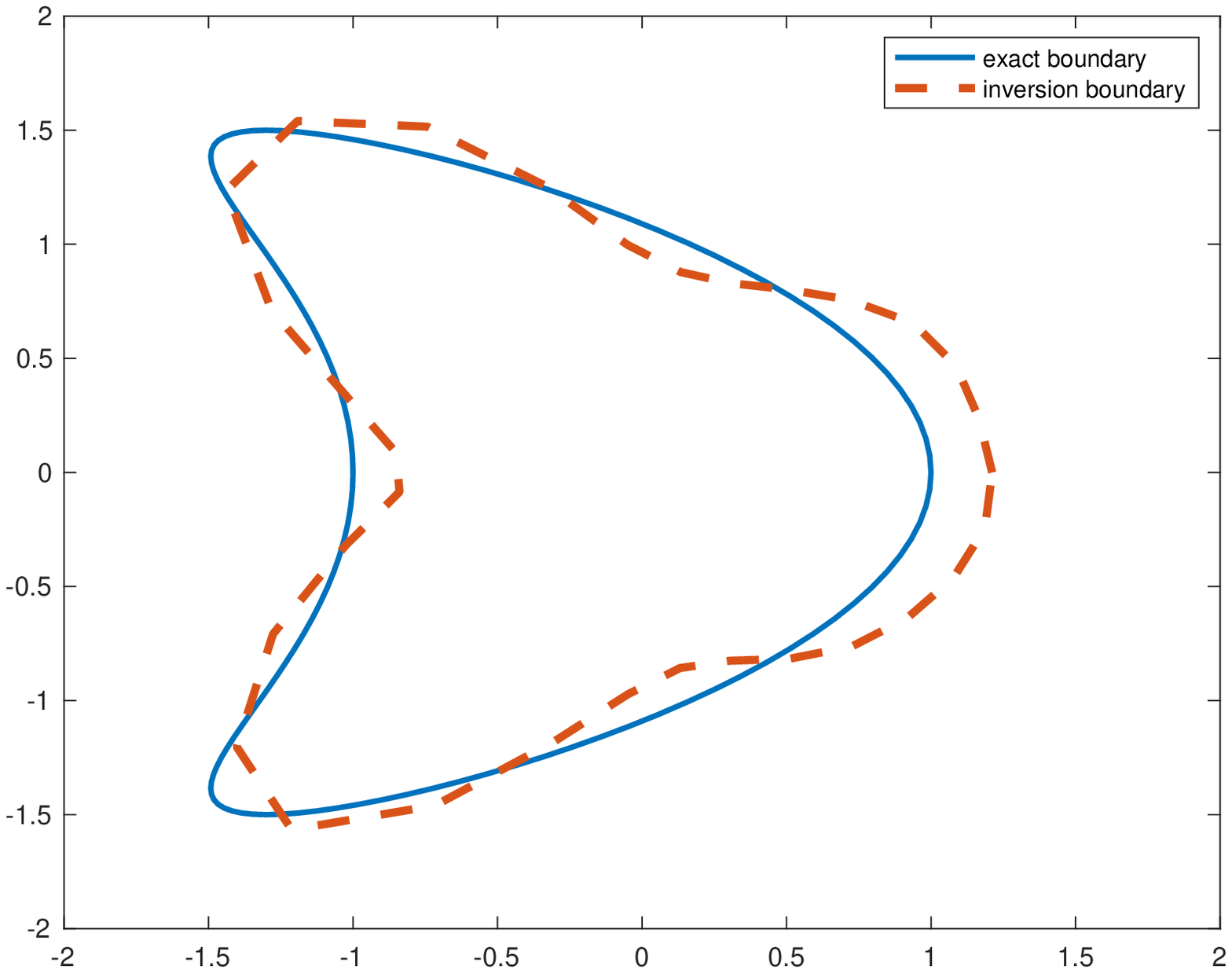}} 
    \subfigure[Roundrect]{\includegraphics[width=4cm]{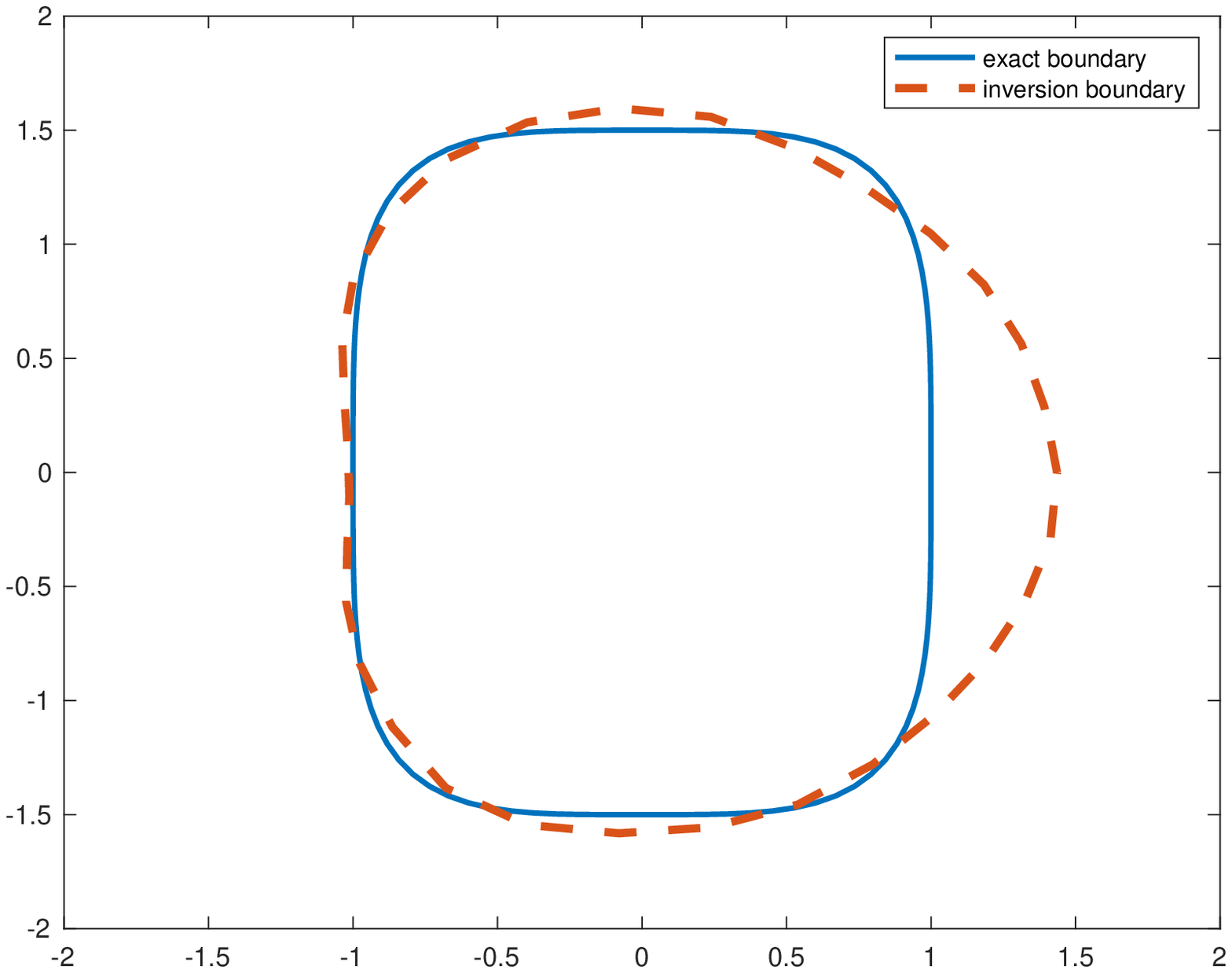}} \\
    \subfigure[Pear]{\includegraphics[width=4cm]{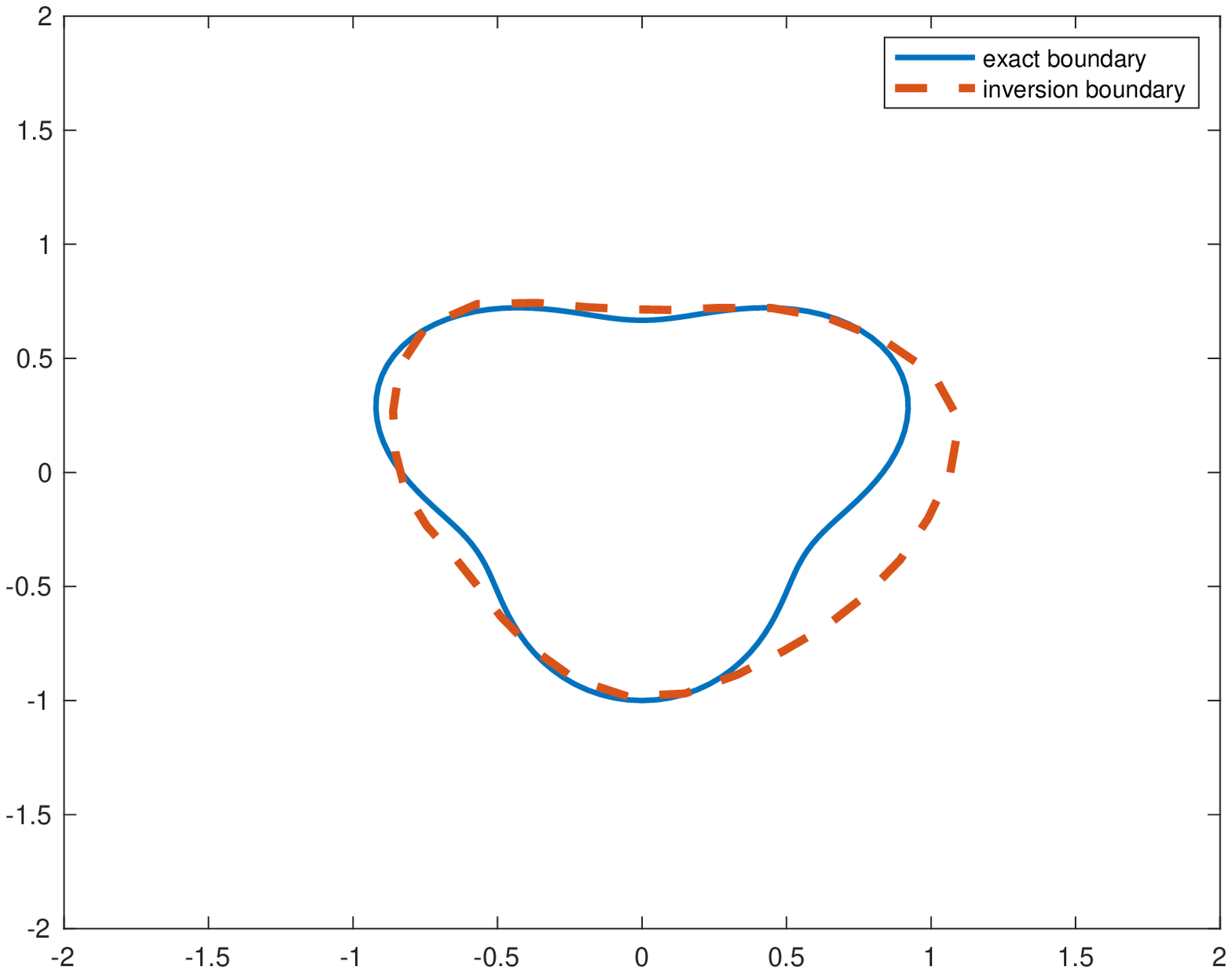}} 
 \subfigure[Acorn]{\includegraphics[width=4cm]{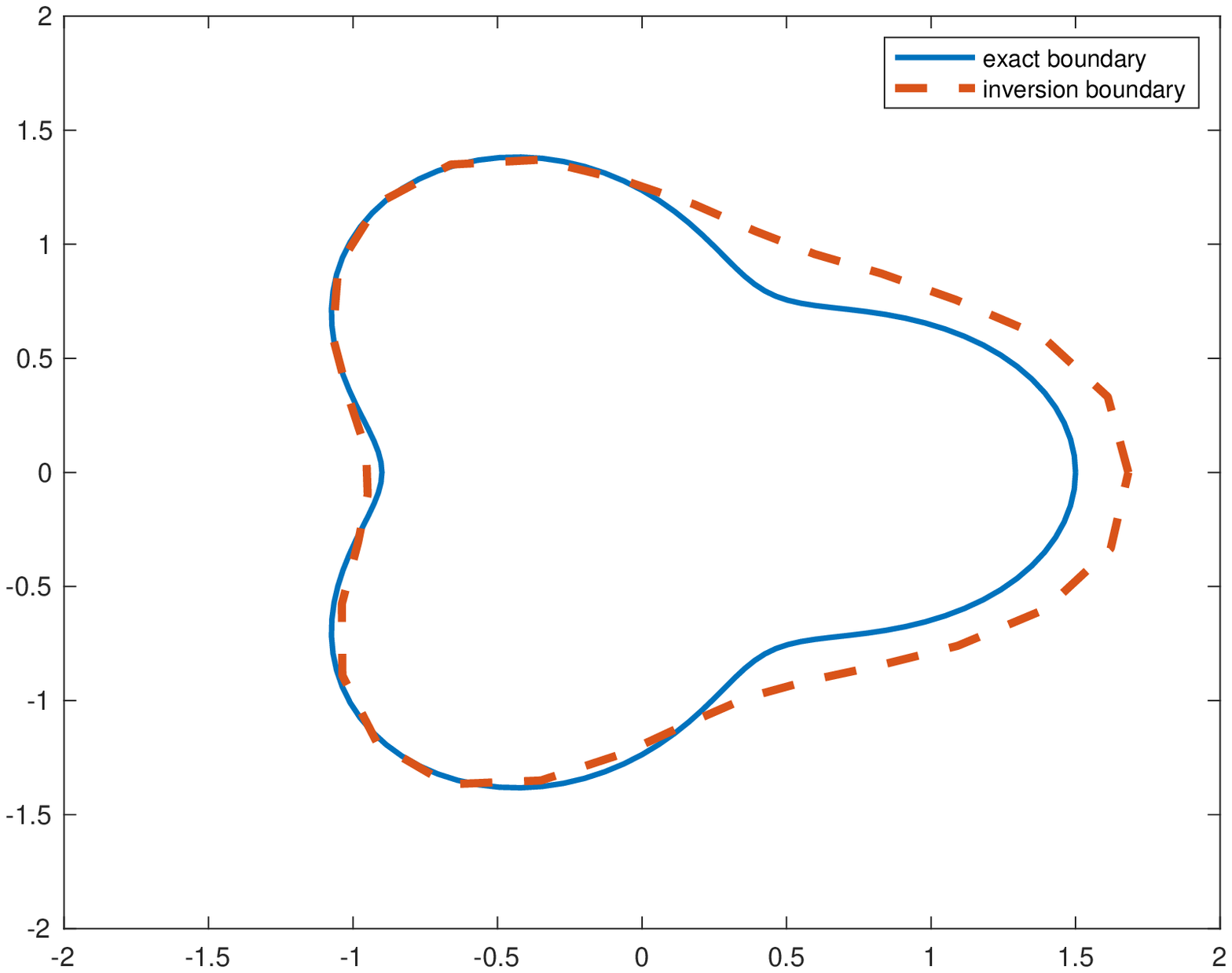}}\\
    \subfigure[Bean]{\includegraphics[width=4cm]{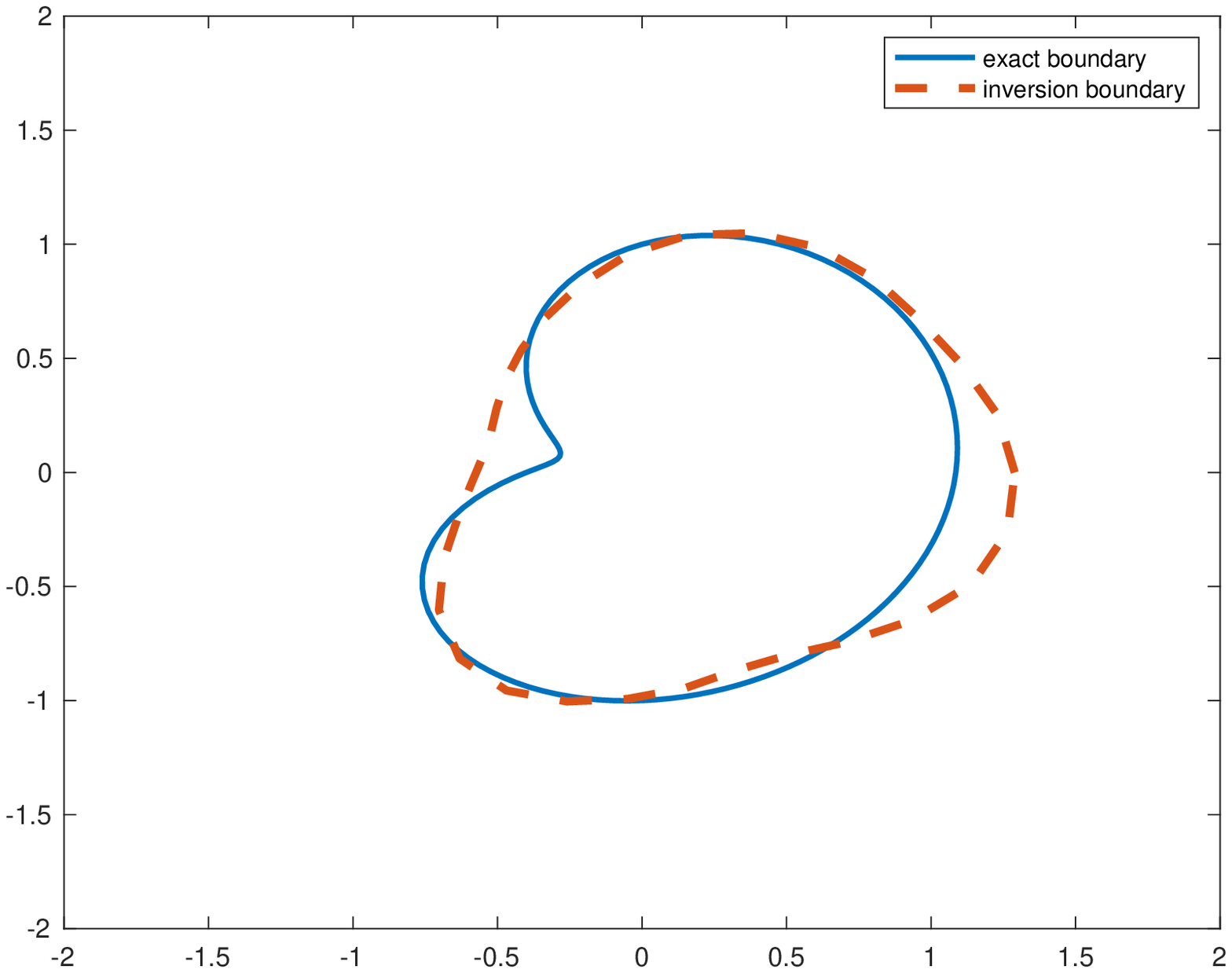}}
     \subfigure[Threelobes]{\includegraphics[width=4cm]{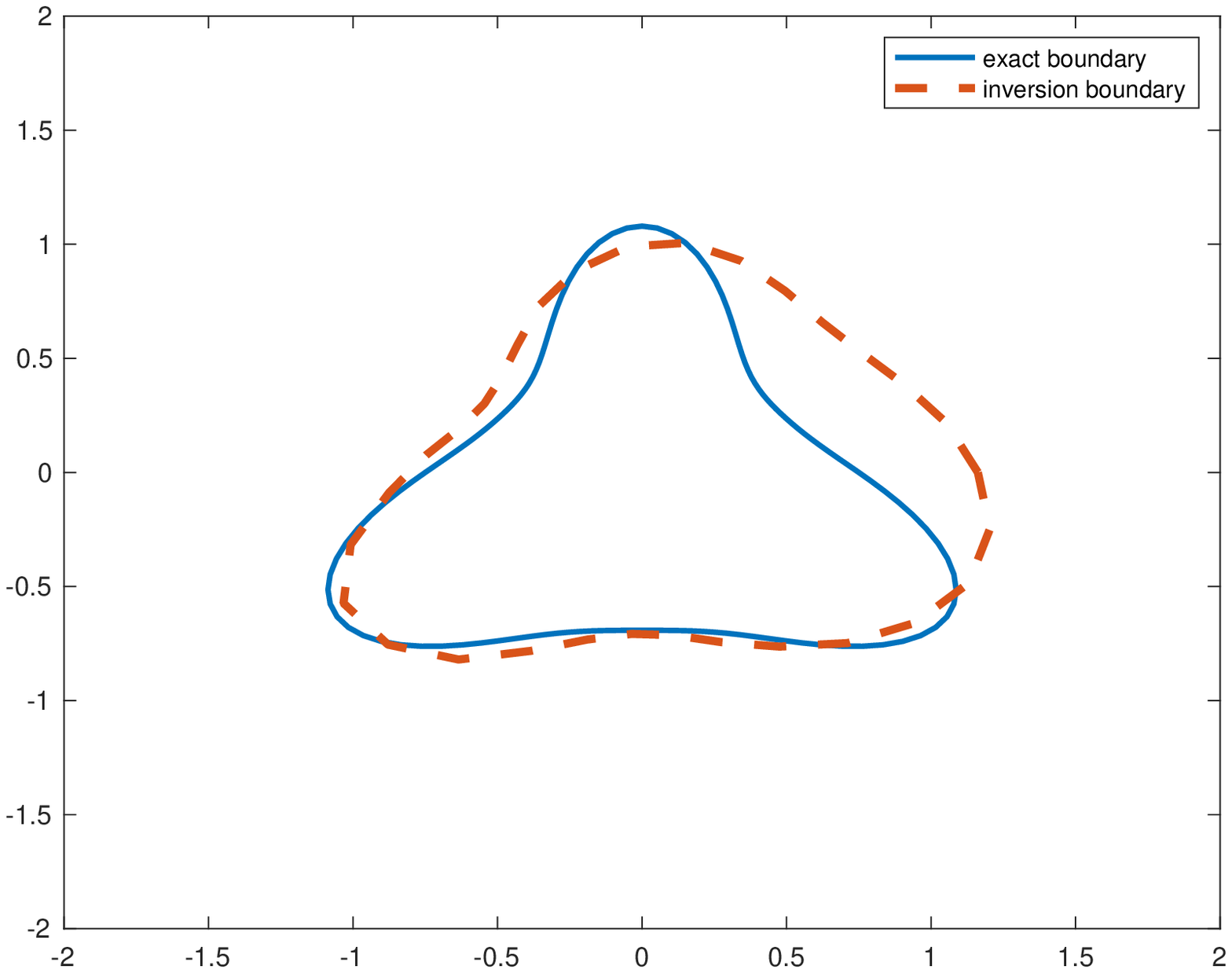}}\\
      \subfigure[Star]{\includegraphics[width=4cm]{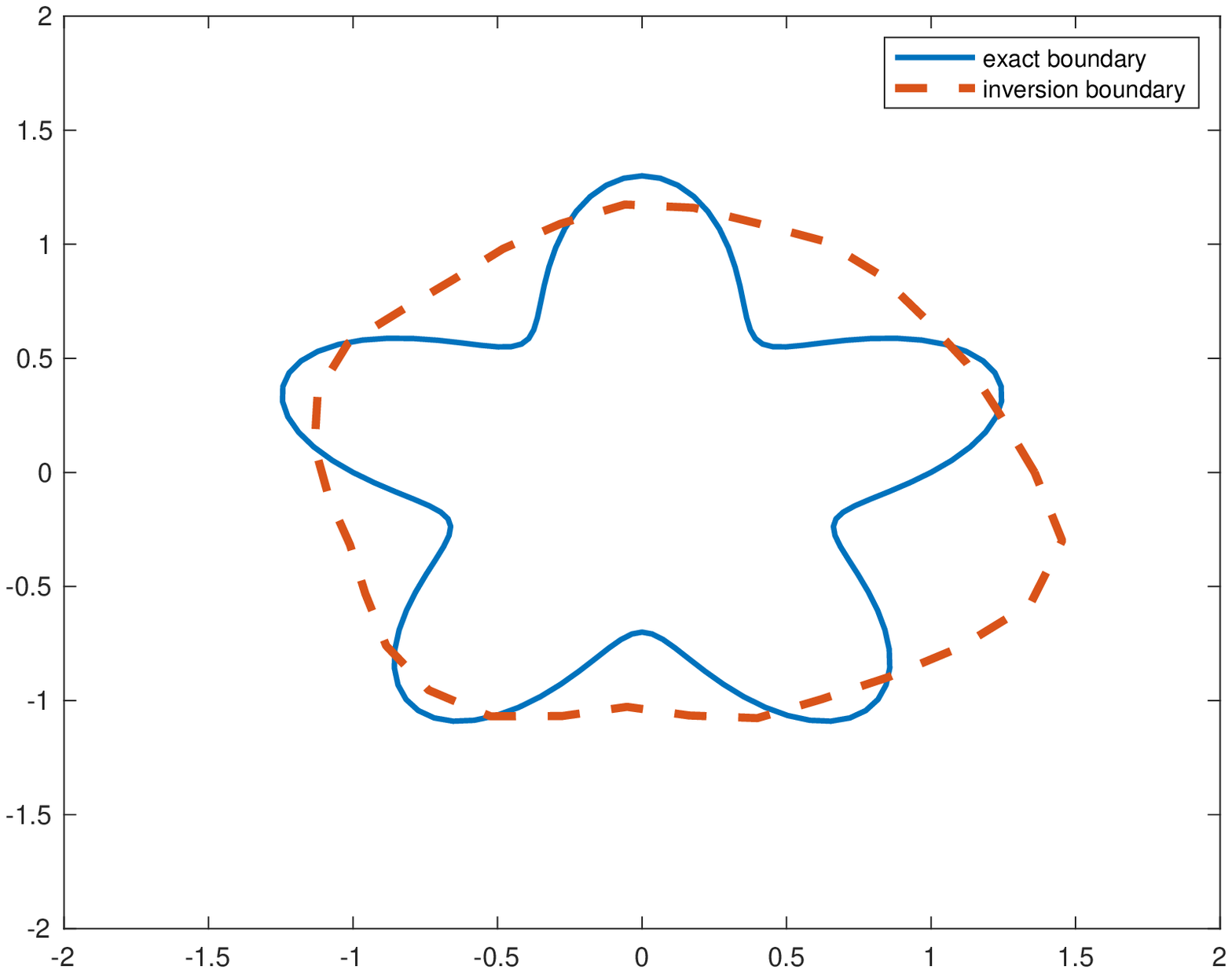}}
       \subfigure[Cloverleaf]{\includegraphics[width=4cm]{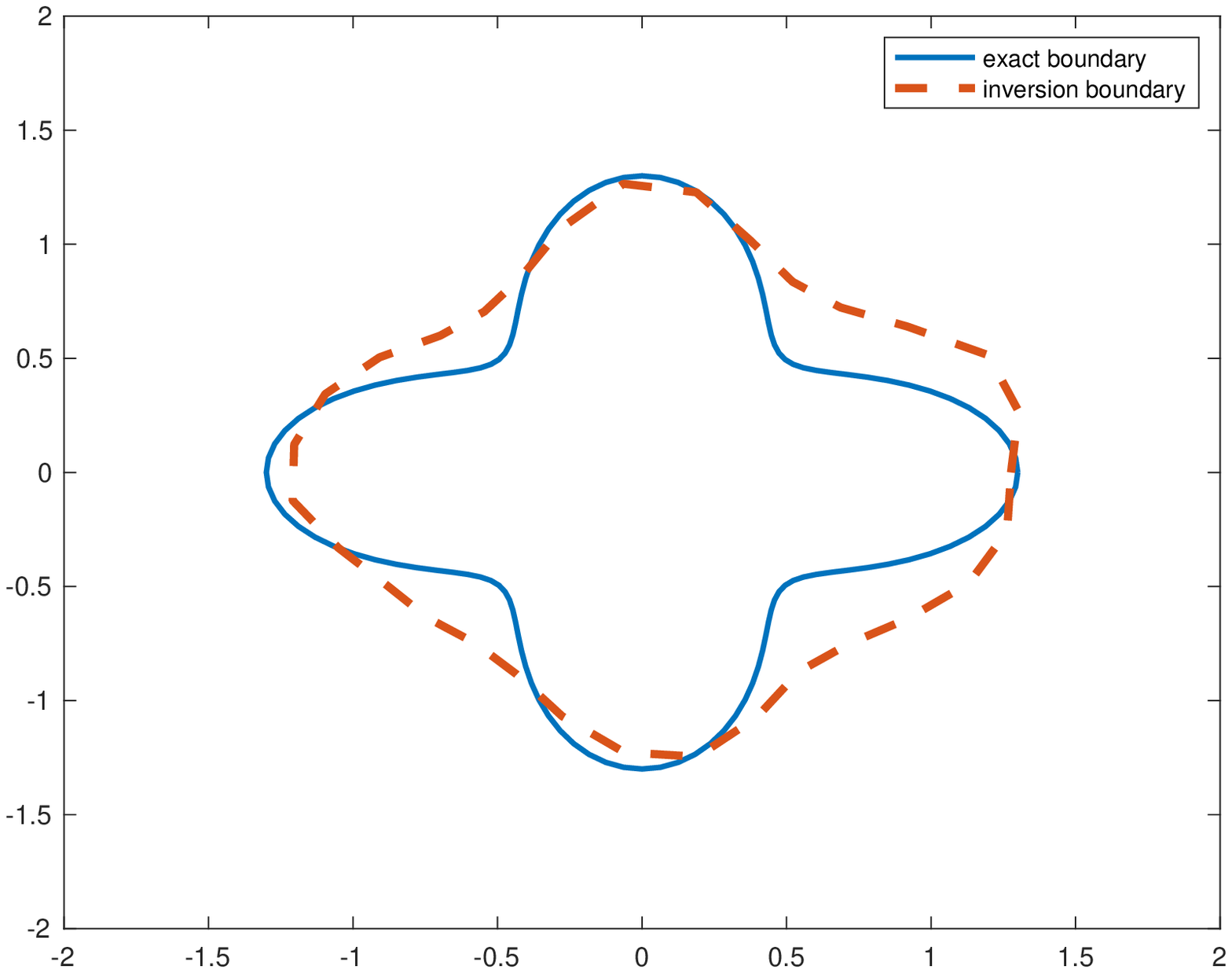}}\\
          \subfigure[Peanut]{\includegraphics[width=4cm]{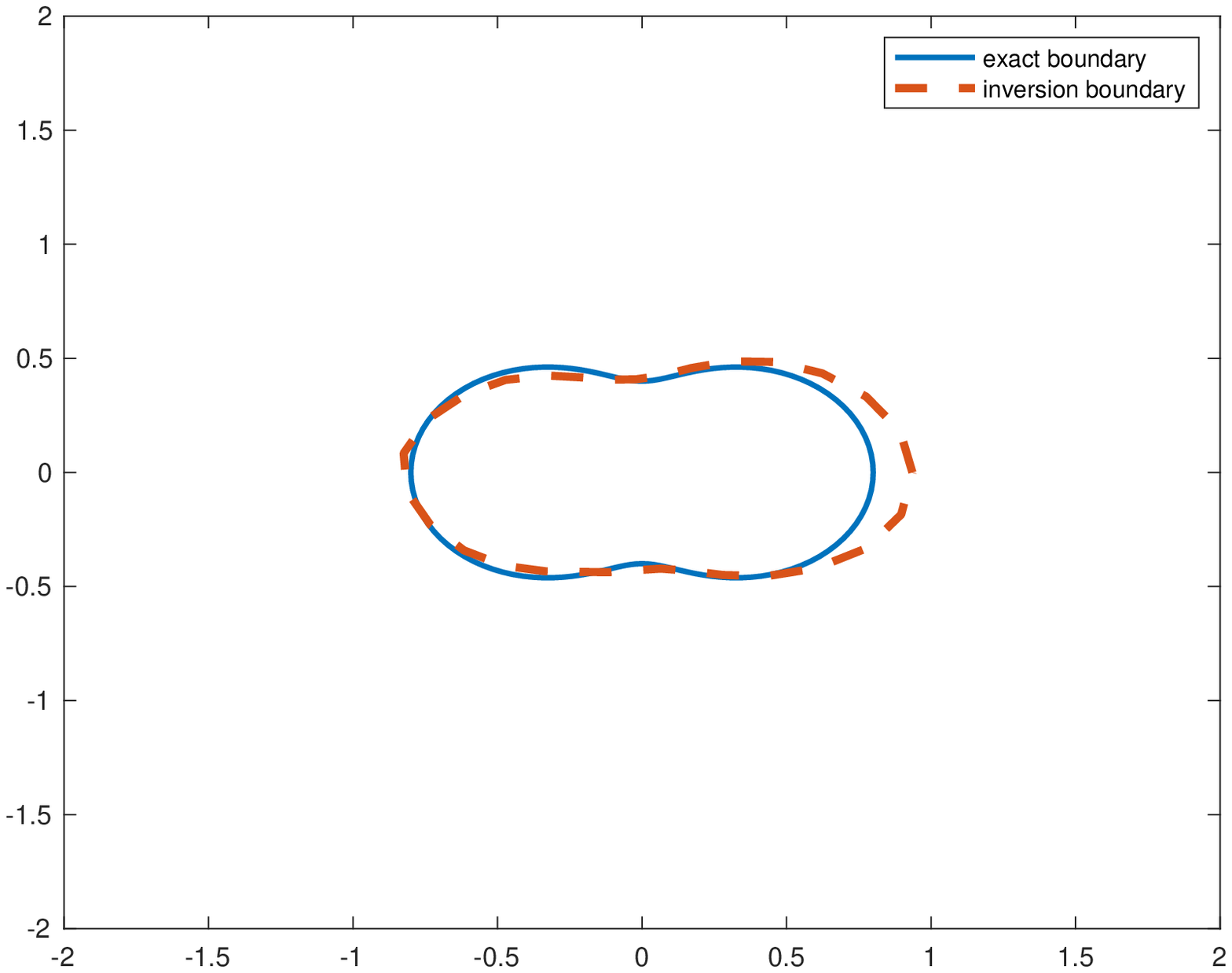}}
       \subfigure[Drop]{\includegraphics[width=4cm]{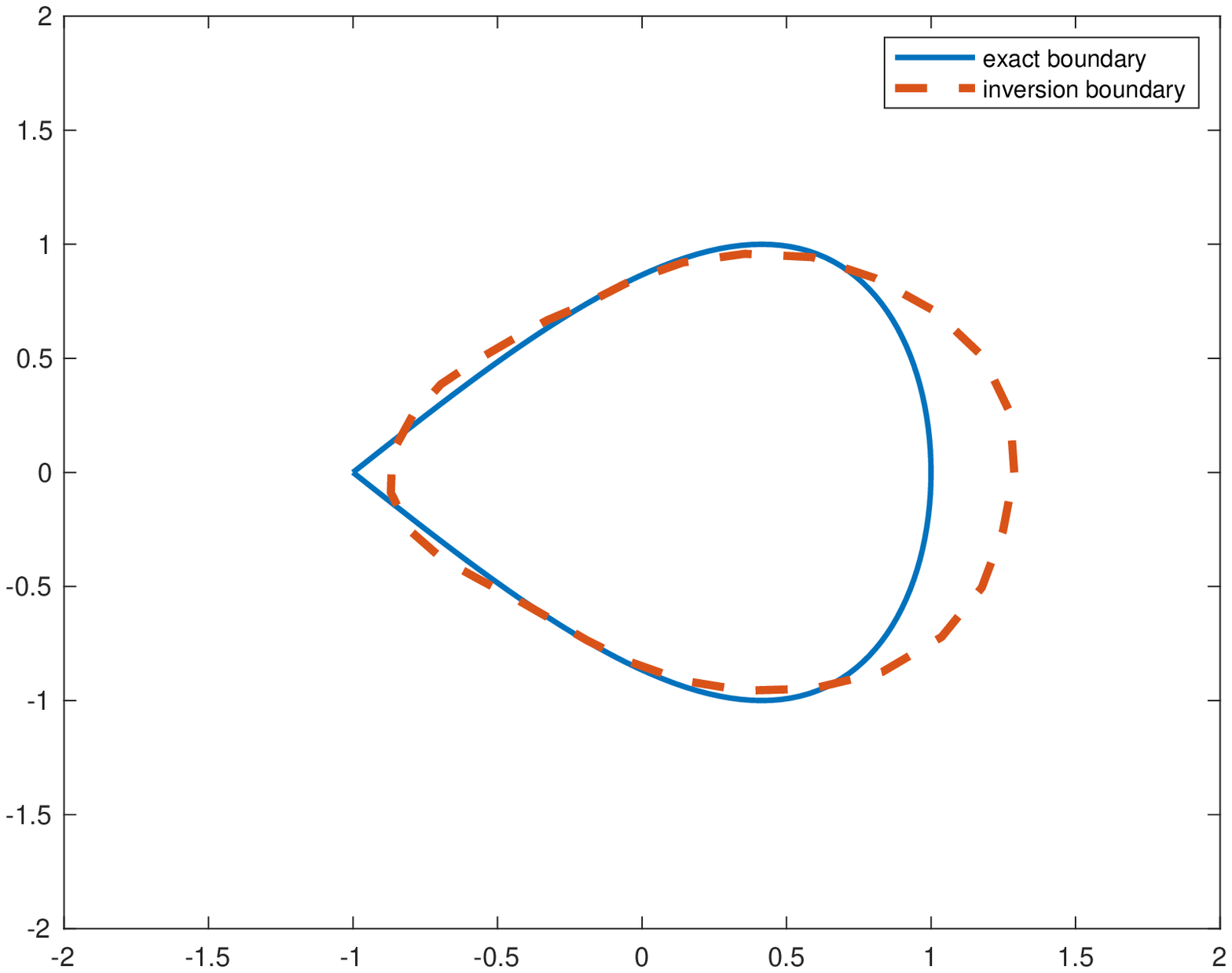}}
   \caption{Reconstructions for different shapes using $u^\infty(\hat{x}, d)$, $(\hat{x}, d)\in \gamma^{\rm o}_{1}\times\gamma^{\rm i}_{1}$ with $\eta_1=1\%, \,\eta_2=1\%$.}
    \label{resul1}
\end{figure}

\begin{figure}
\captionsetup[subfigure]{labelformat=empty}
    \centering
    \subfigure[Kite]{\includegraphics[width=4cm]{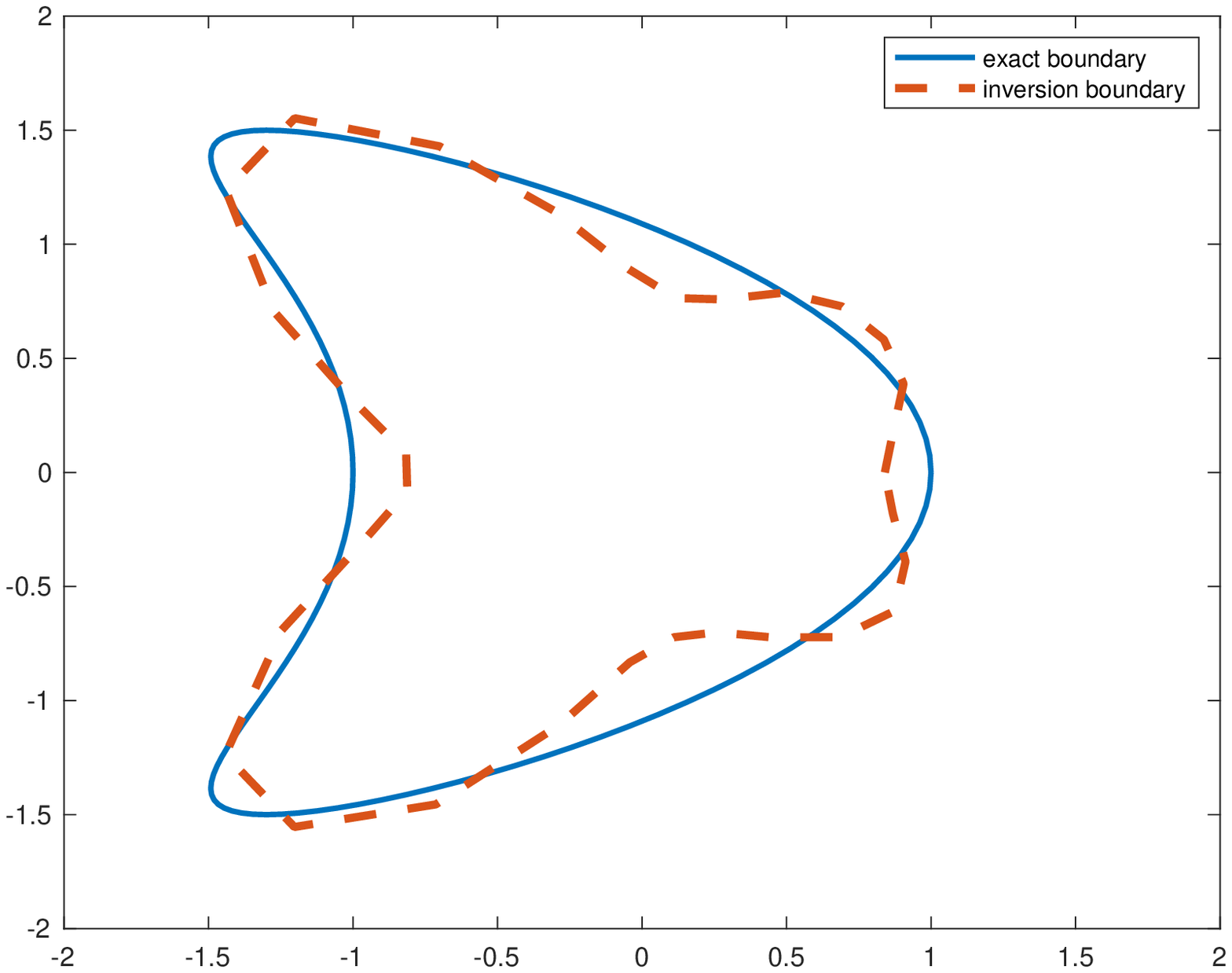}} 
    \subfigure[Roundrect]{\includegraphics[width=4cm]{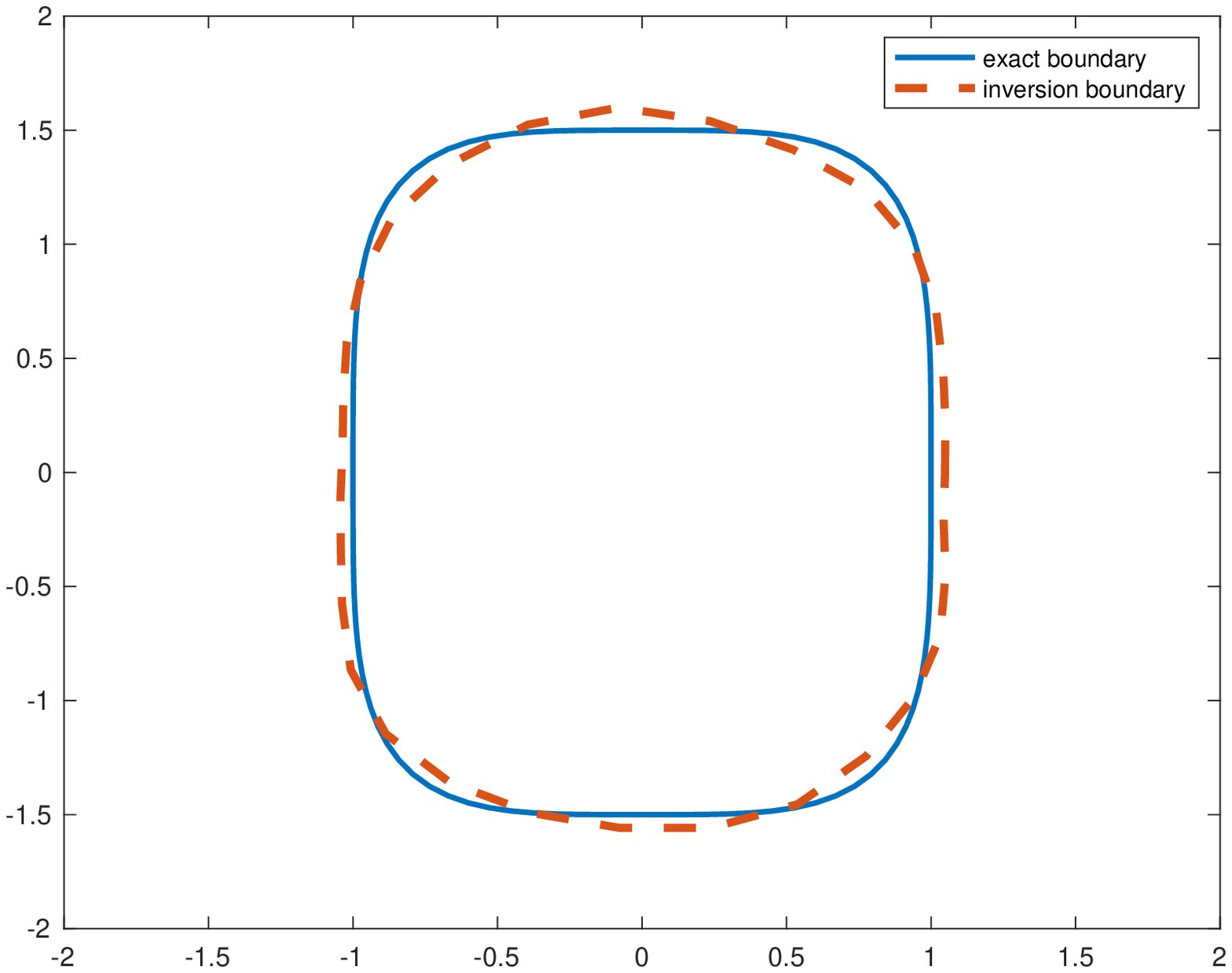}} \\
    \subfigure[Pear]{\includegraphics[width=4cm]{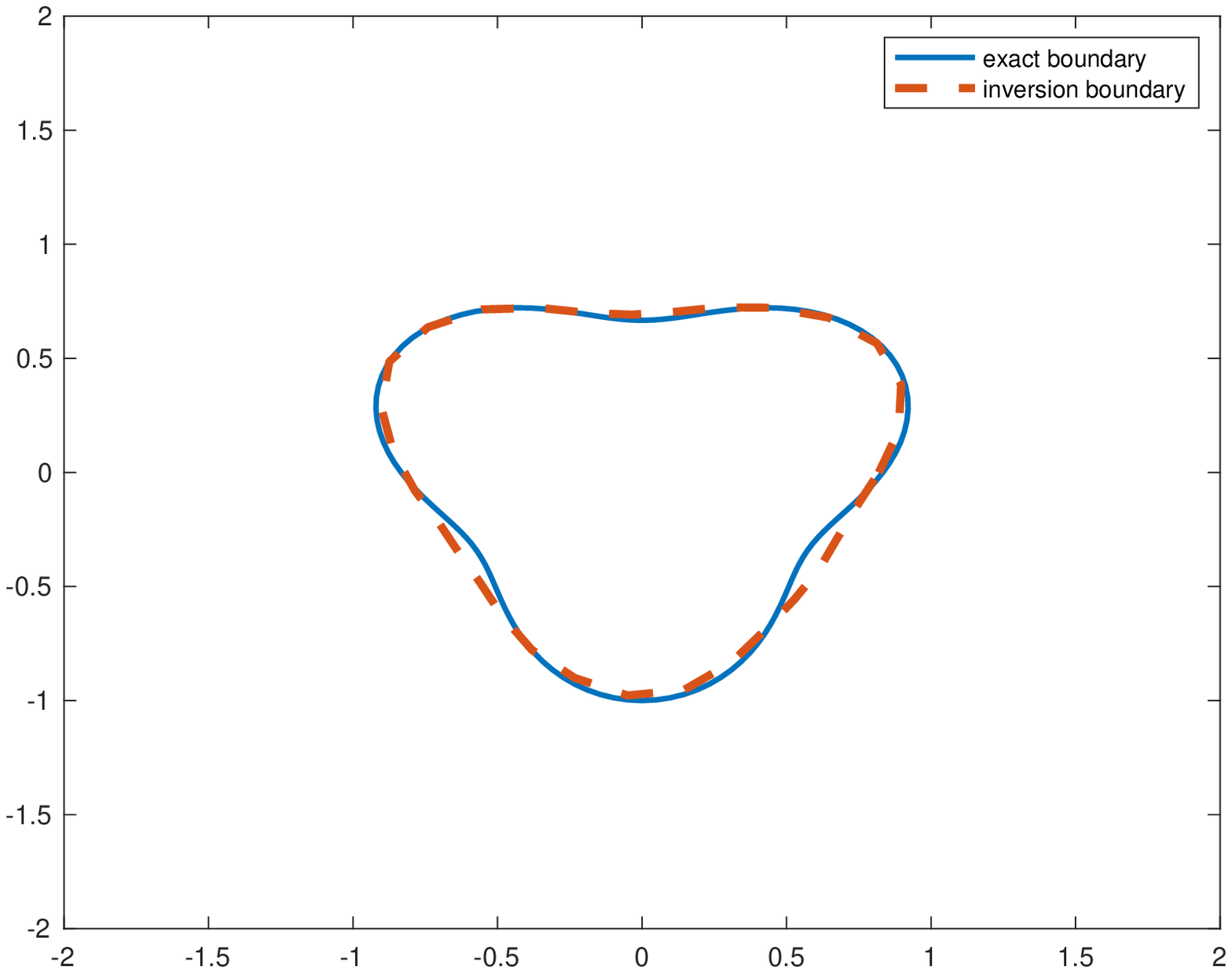}} 
 \subfigure[Acorn]{\includegraphics[width=4cm]{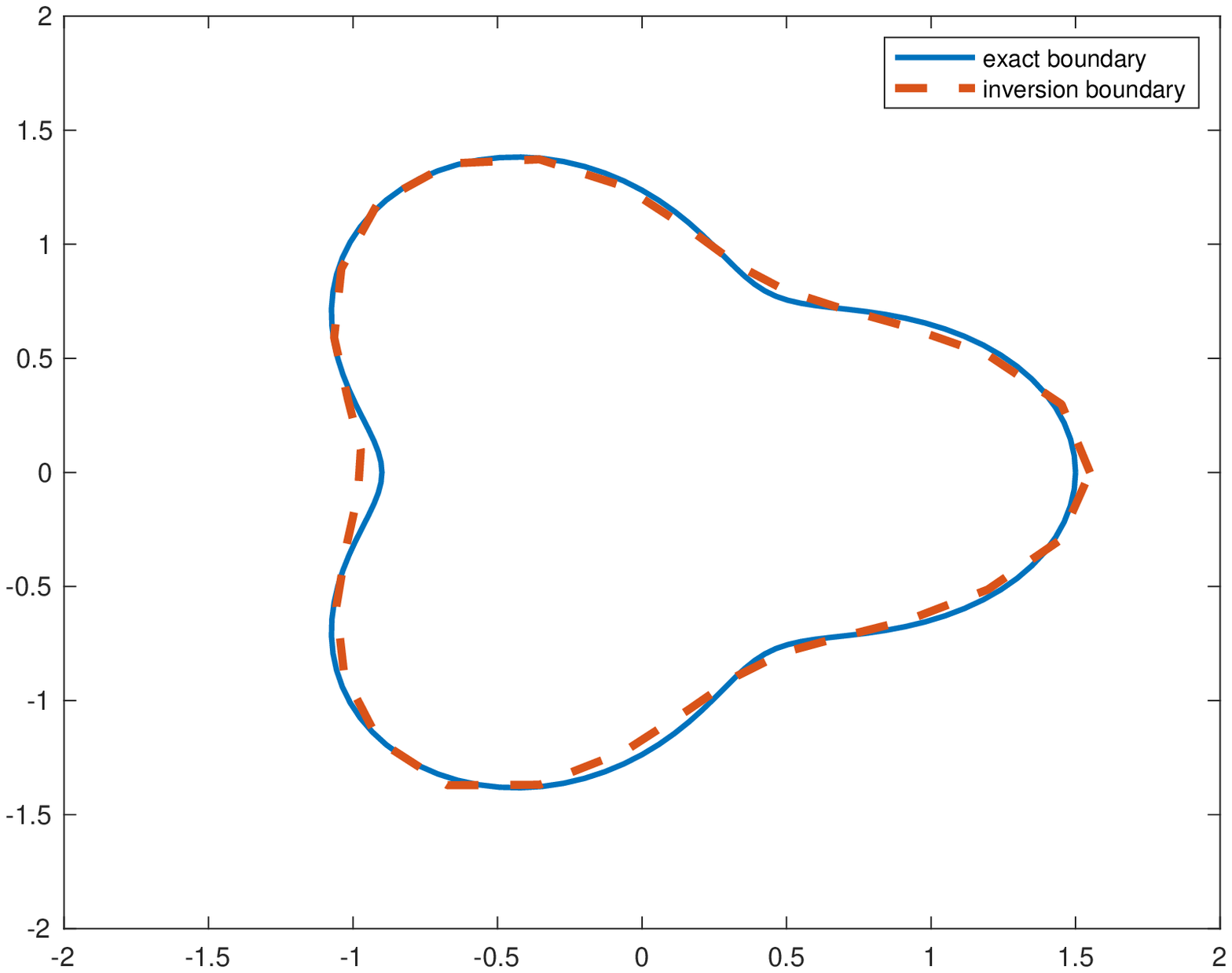}}\\
    \subfigure[Bean]{\includegraphics[width=4cm]{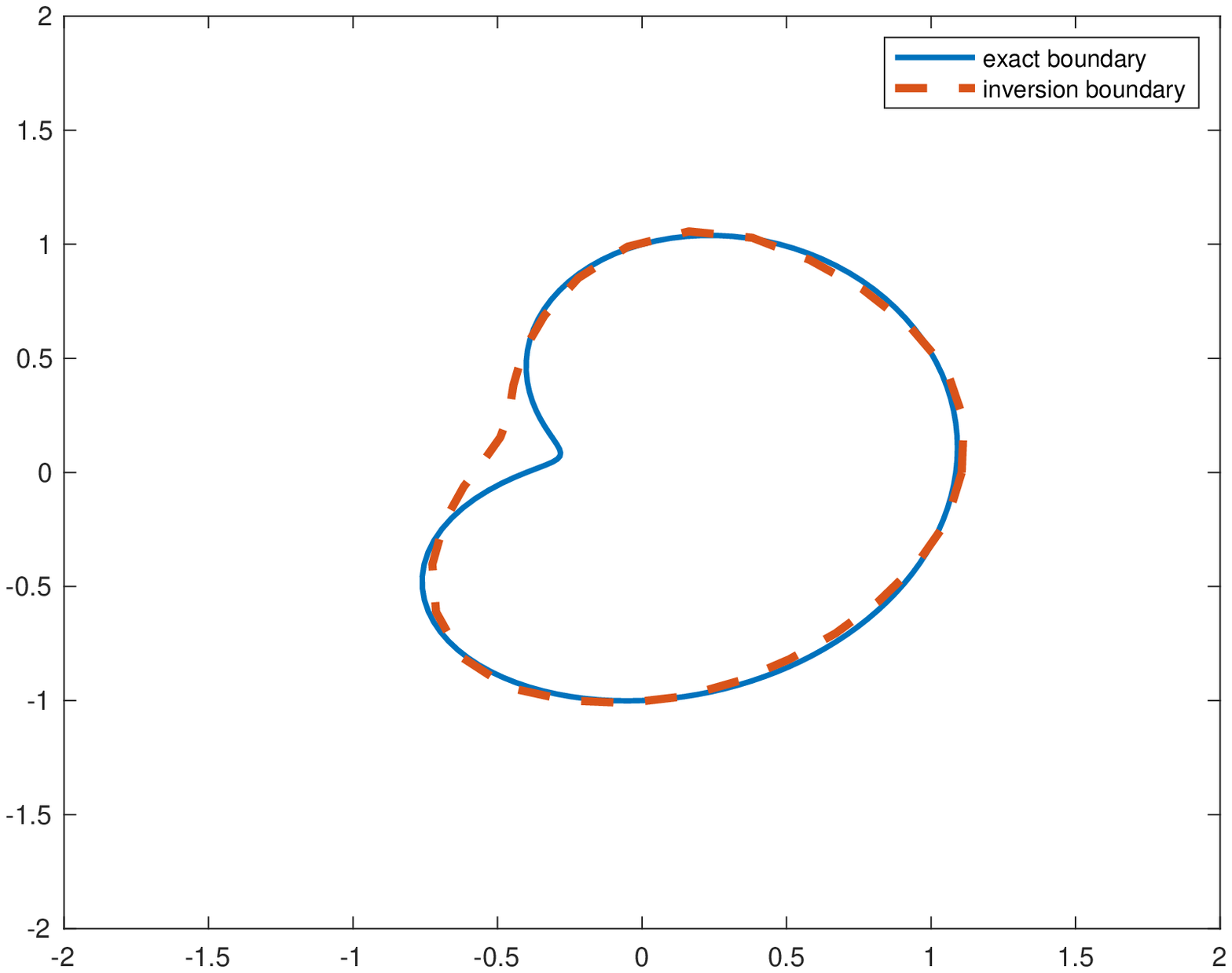}}
     \subfigure[Threelobes]{\includegraphics[width=4cm]{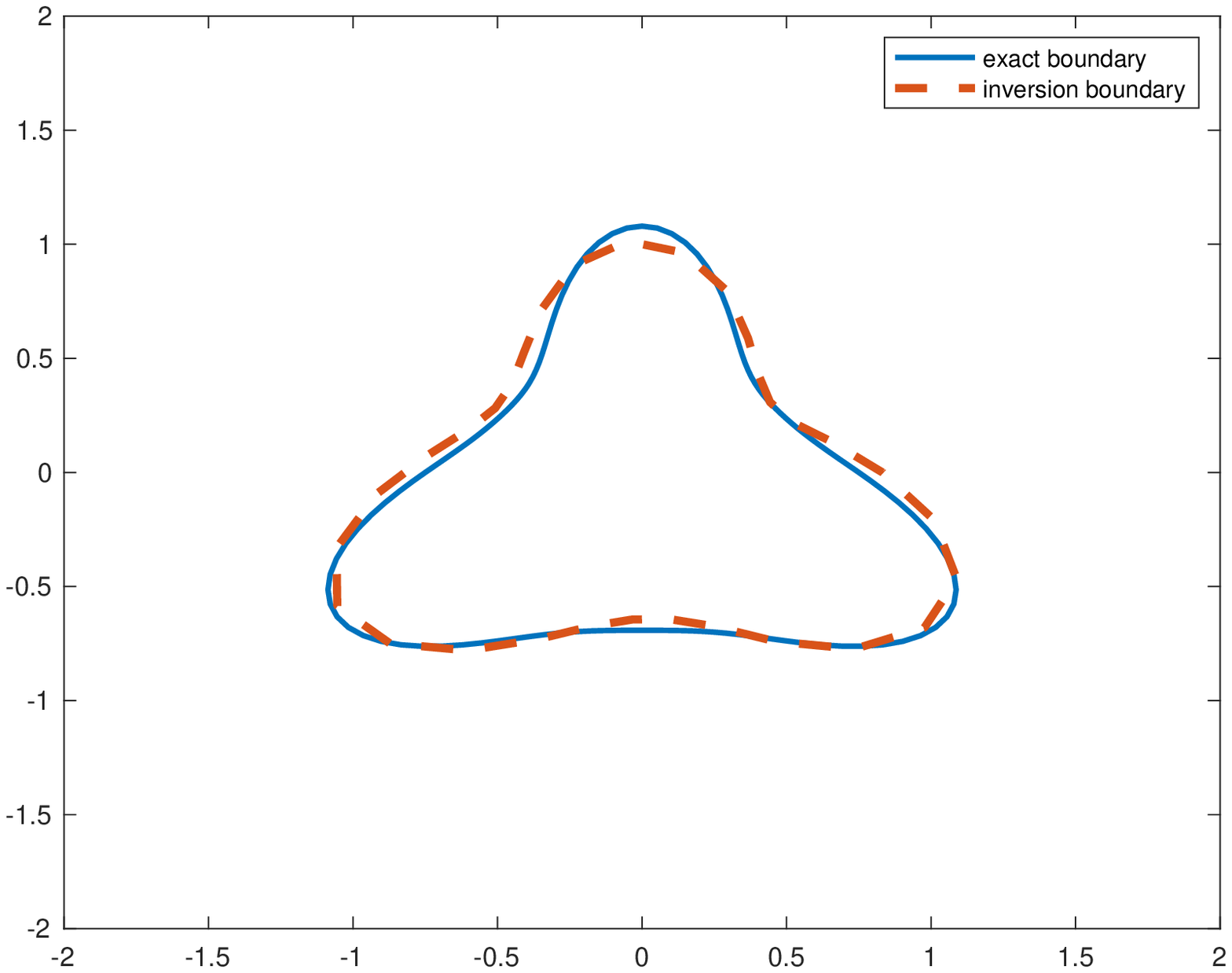}}\\
      \subfigure[Star]{\includegraphics[width=4cm]{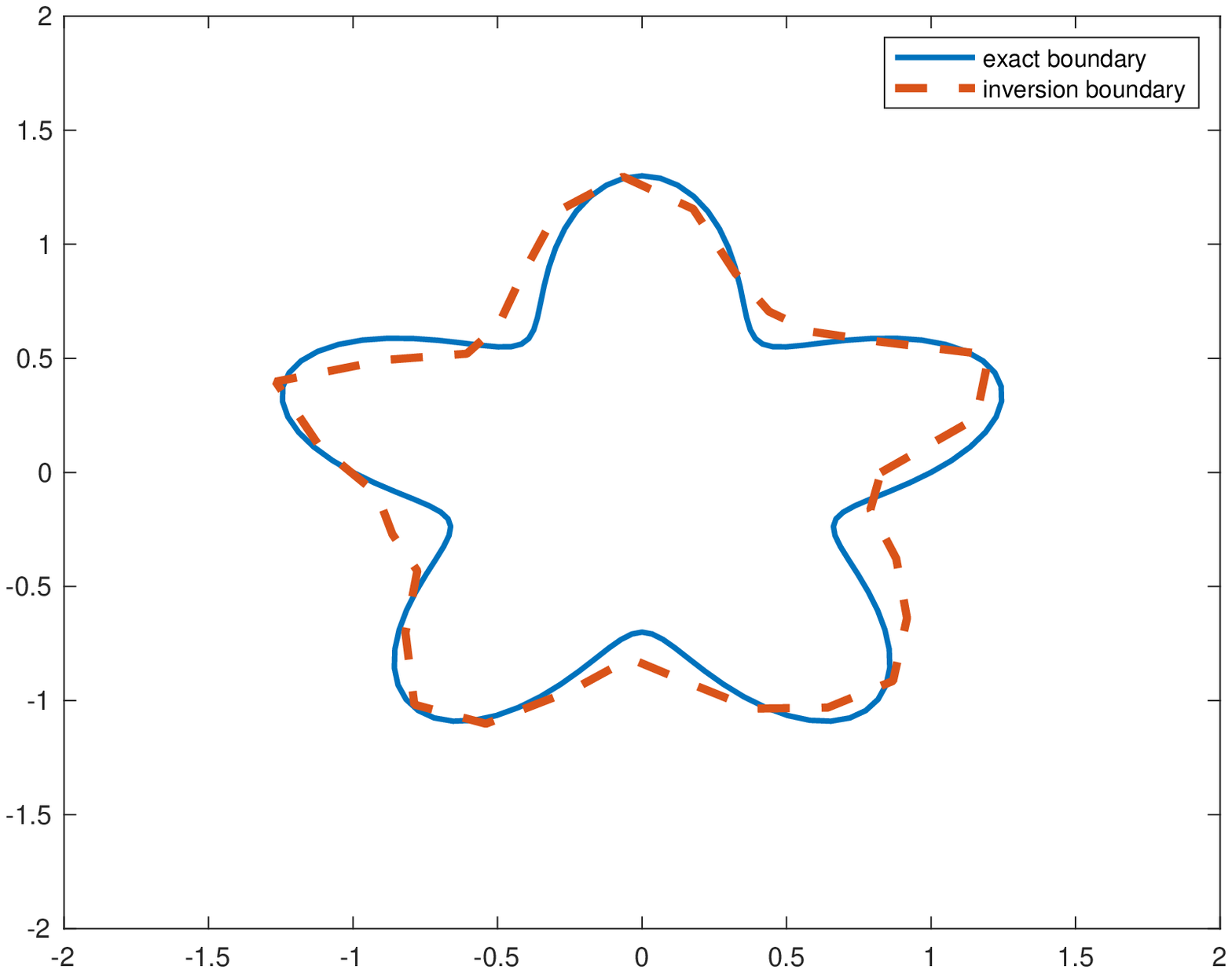}}
       \subfigure[Cloverleaf]{\includegraphics[width=4cm]{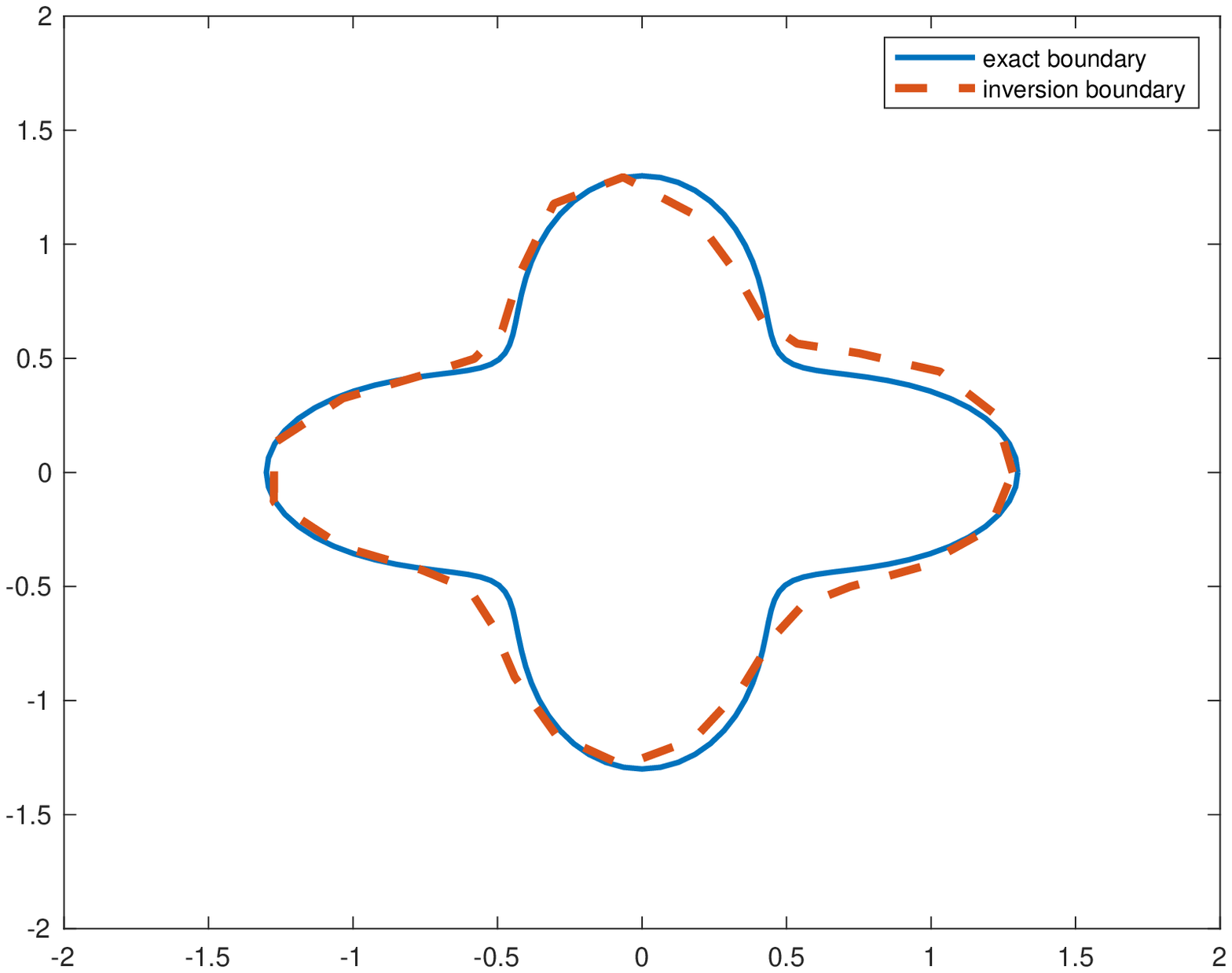}}\\
          \subfigure[Peanut]{\includegraphics[width=4cm]{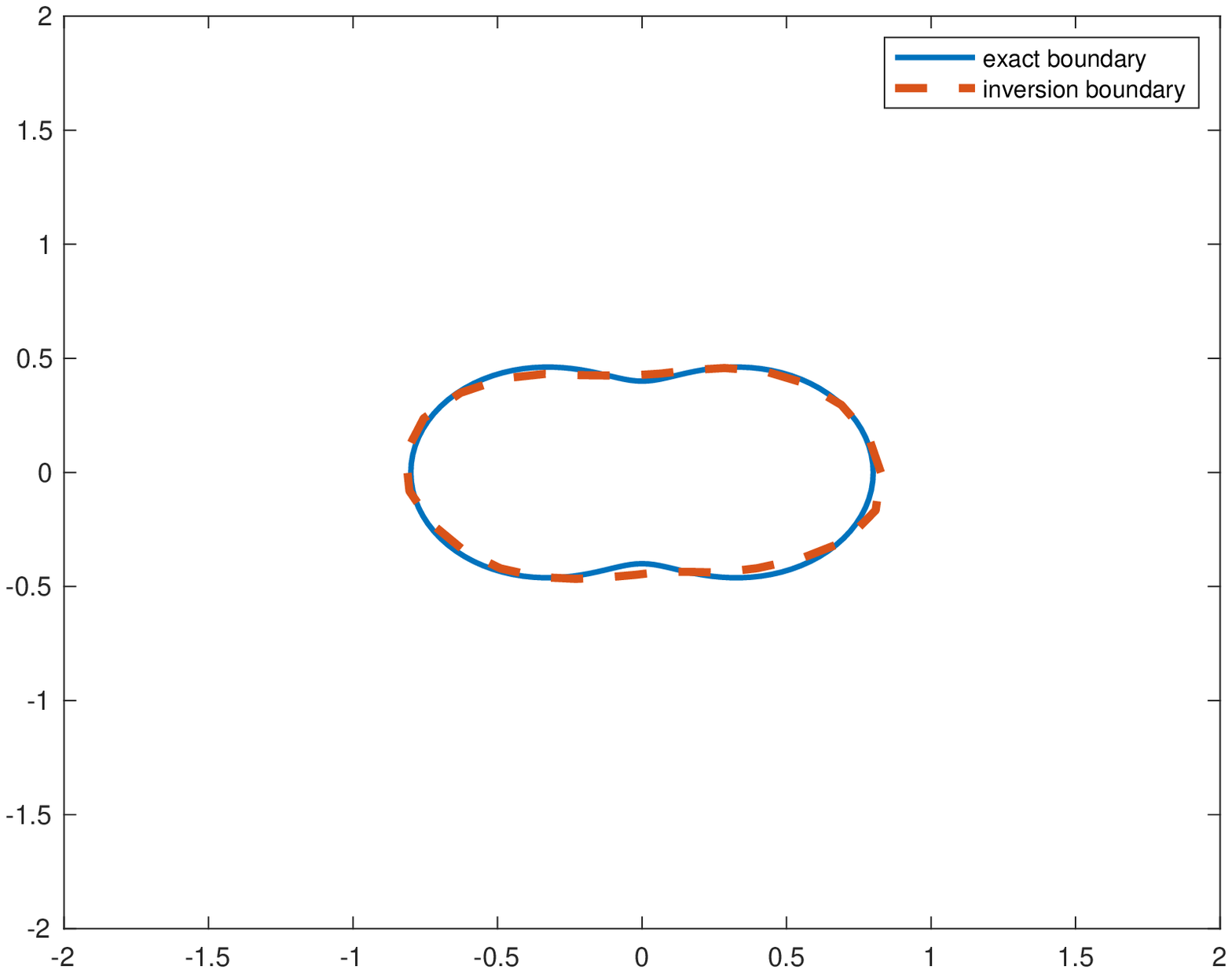}}
       \subfigure[Drop]{\includegraphics[width=4cm]{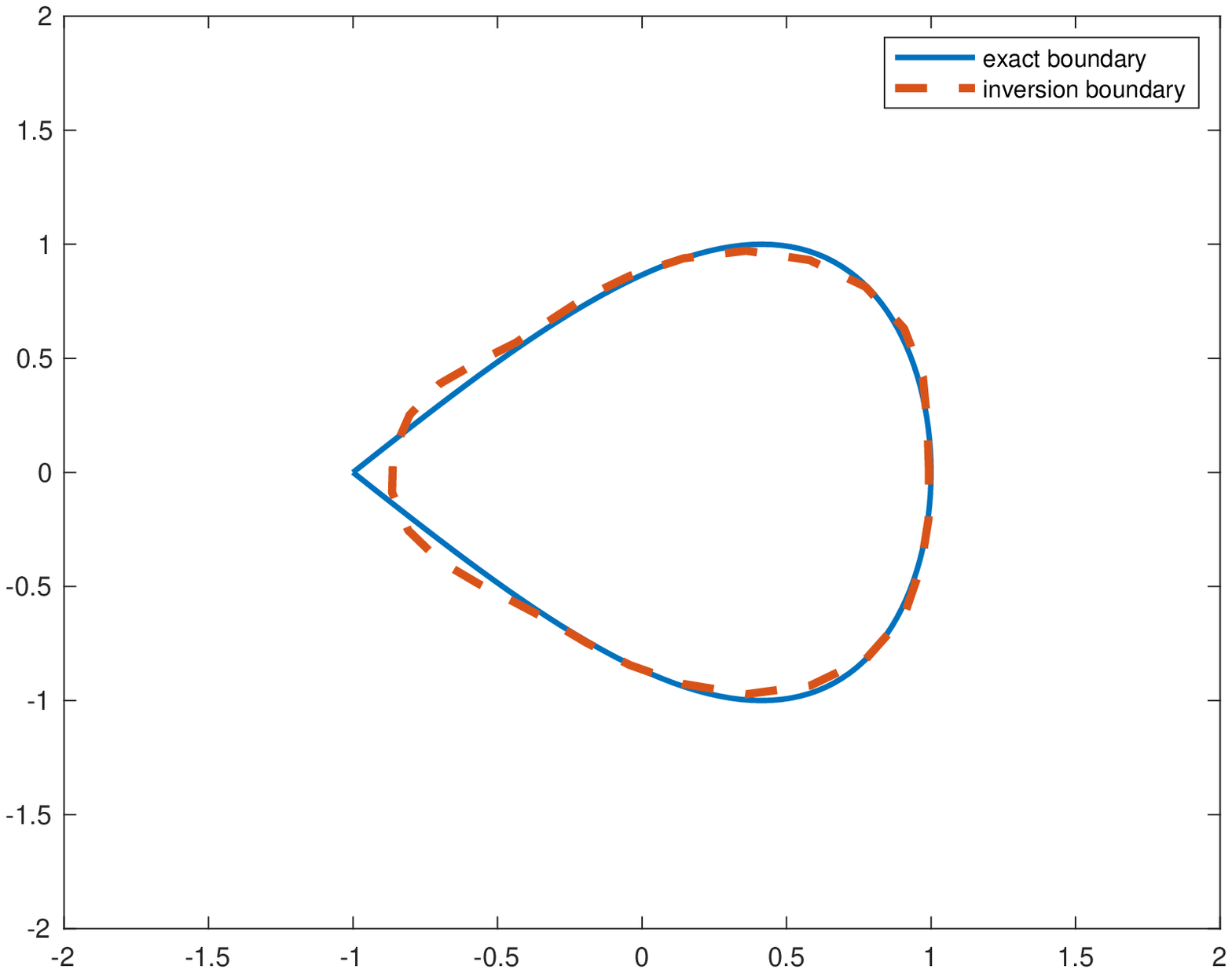}}
   \caption{Reconstructions for different shapes using $u^\infty(\hat{x}, d)$, $(\hat{x}, d)\in \gamma^{\rm o}_{1}\times\gamma^{\rm i}_{2}$ with $\eta_1=1\%, \,\eta_2=1\%$.}
    \label{resul2}
\end{figure}

We first consider the cases of one incident wave and two incident waves, i.e., $\gamma^{\rm i}=\gamma_1^{\rm i}$ and $\gamma^{\rm i}=\gamma_2^{\rm i}$, with full observation aperture. Taking one observation aperture with $\gamma_{1}^{\rm o}$, we display the numerical reconstructions in Fig. \ref{resul1} ($\gamma^{\rm o}\times \gamma^{\rm i}=\gamma_1^{\rm o}\times \gamma_{1}^{\rm i}$), Fig. \ref{resul2} ($\gamma^{\rm o}\times\gamma^{\rm i}=\gamma_1^{\rm o}\times\gamma_{2}^{\rm i}$)  using noise level $\eta_1=1\%, \eta_2=1\%$, respectively. 
The results show that the two incident waves perform better. In the following tests, we use two incident waves case. 

\begin{figure}
\captionsetup[subfigure]{labelformat=empty}
    \centering
    \subfigure[Kite]{\includegraphics[width=4cm]{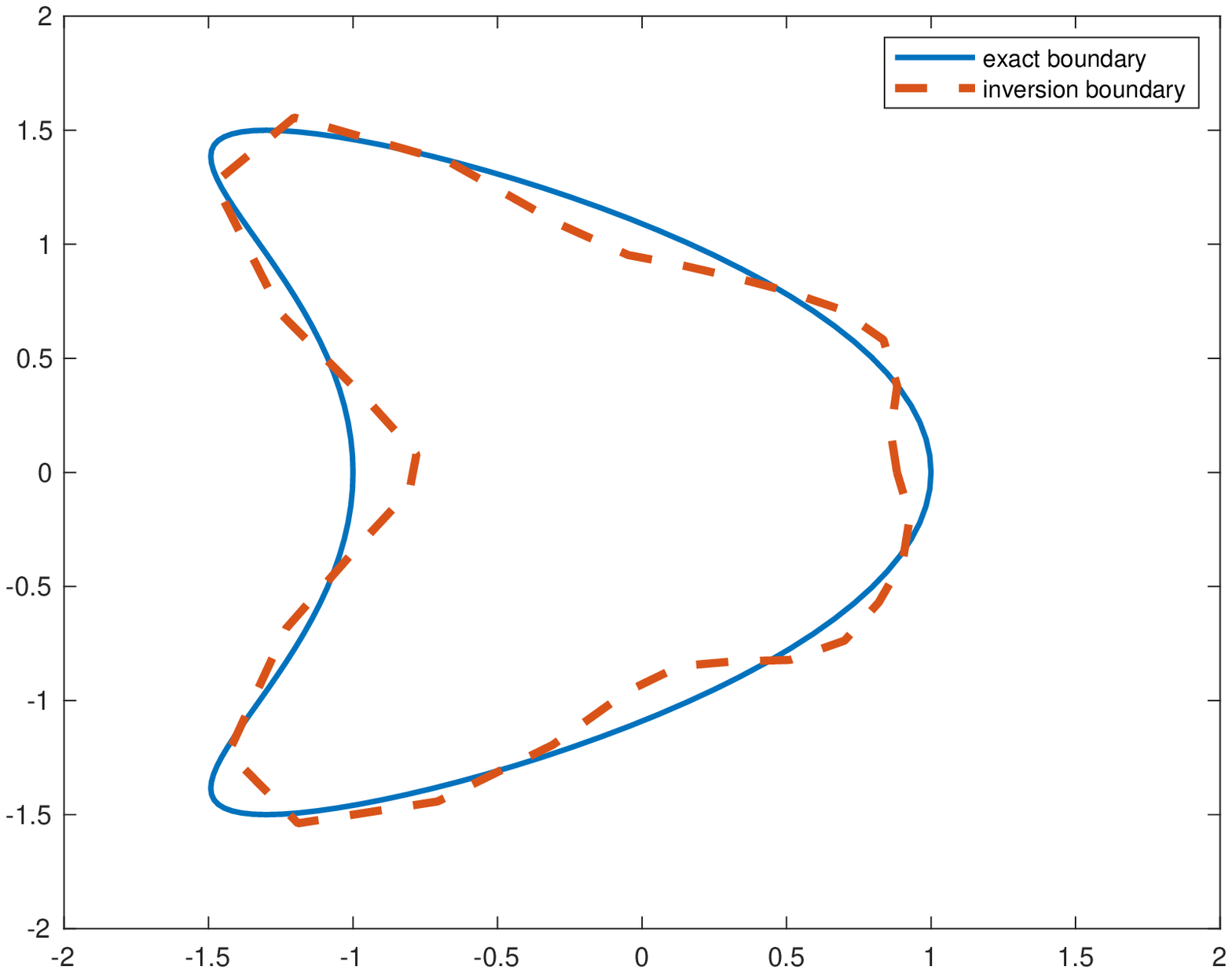}} 
    \subfigure[Roundrect]{\includegraphics[width=4cm]{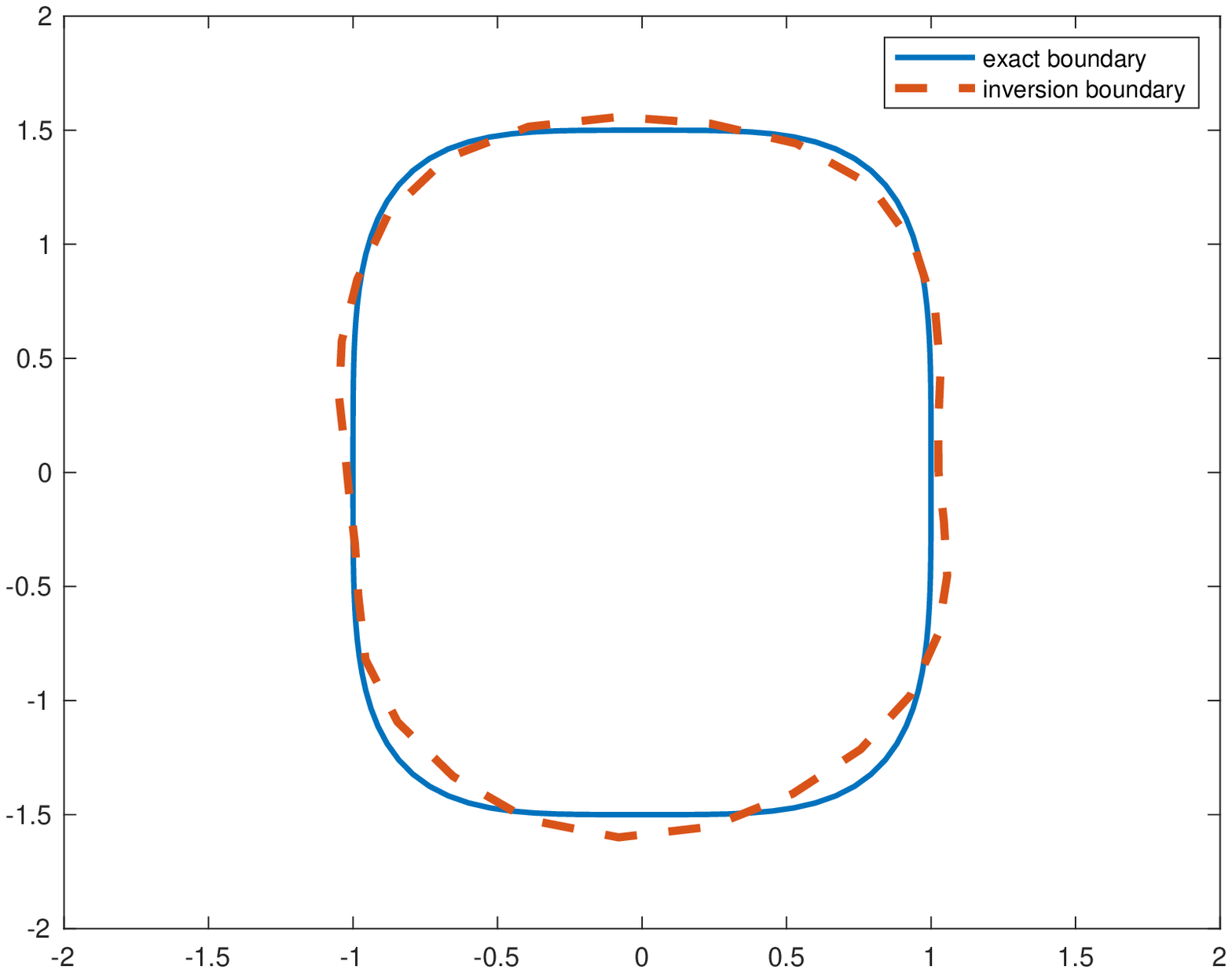}} \\
    \subfigure[Pear]{\includegraphics[width=4cm]{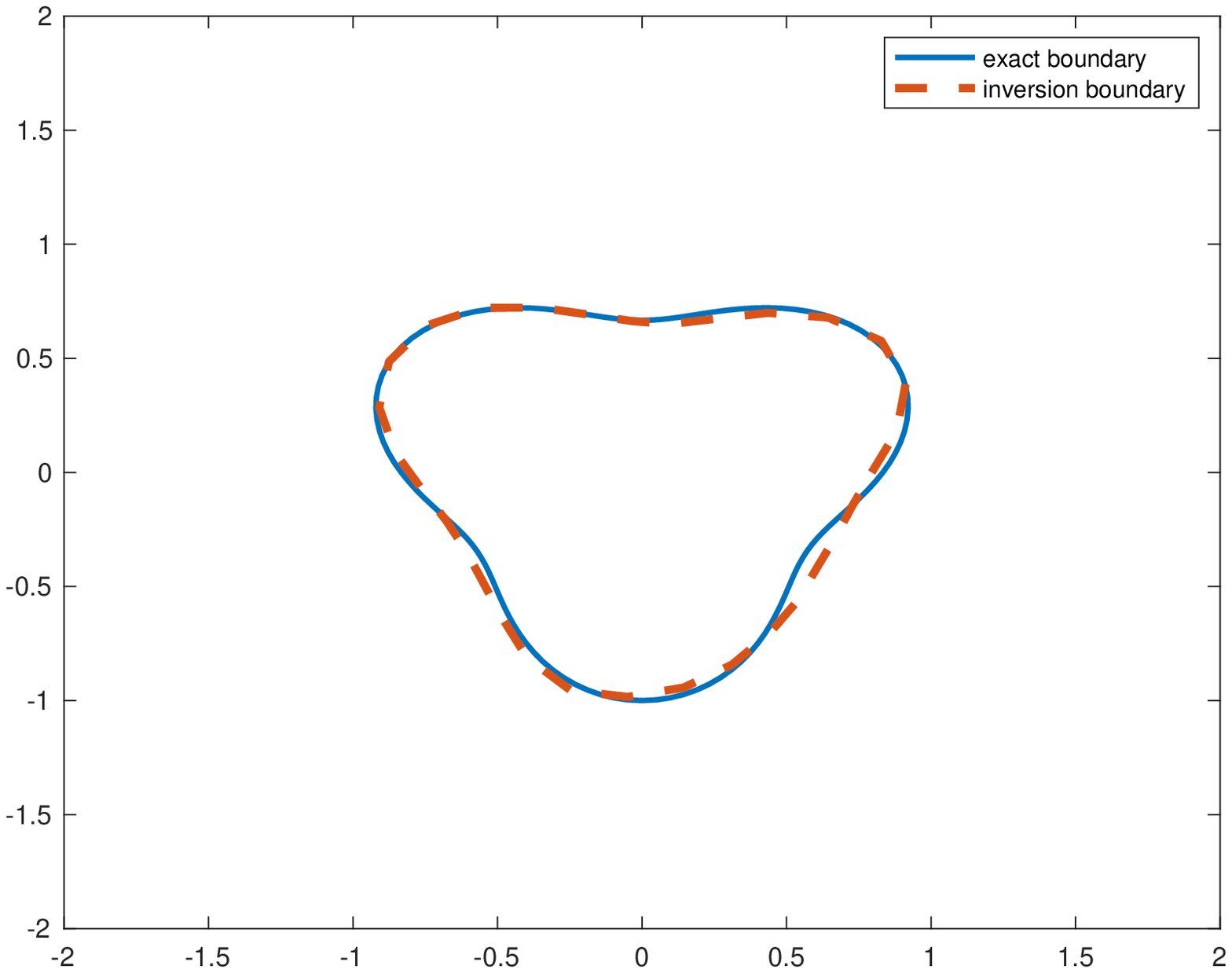}} 
 \subfigure[Acorn]{\includegraphics[width=4cm]{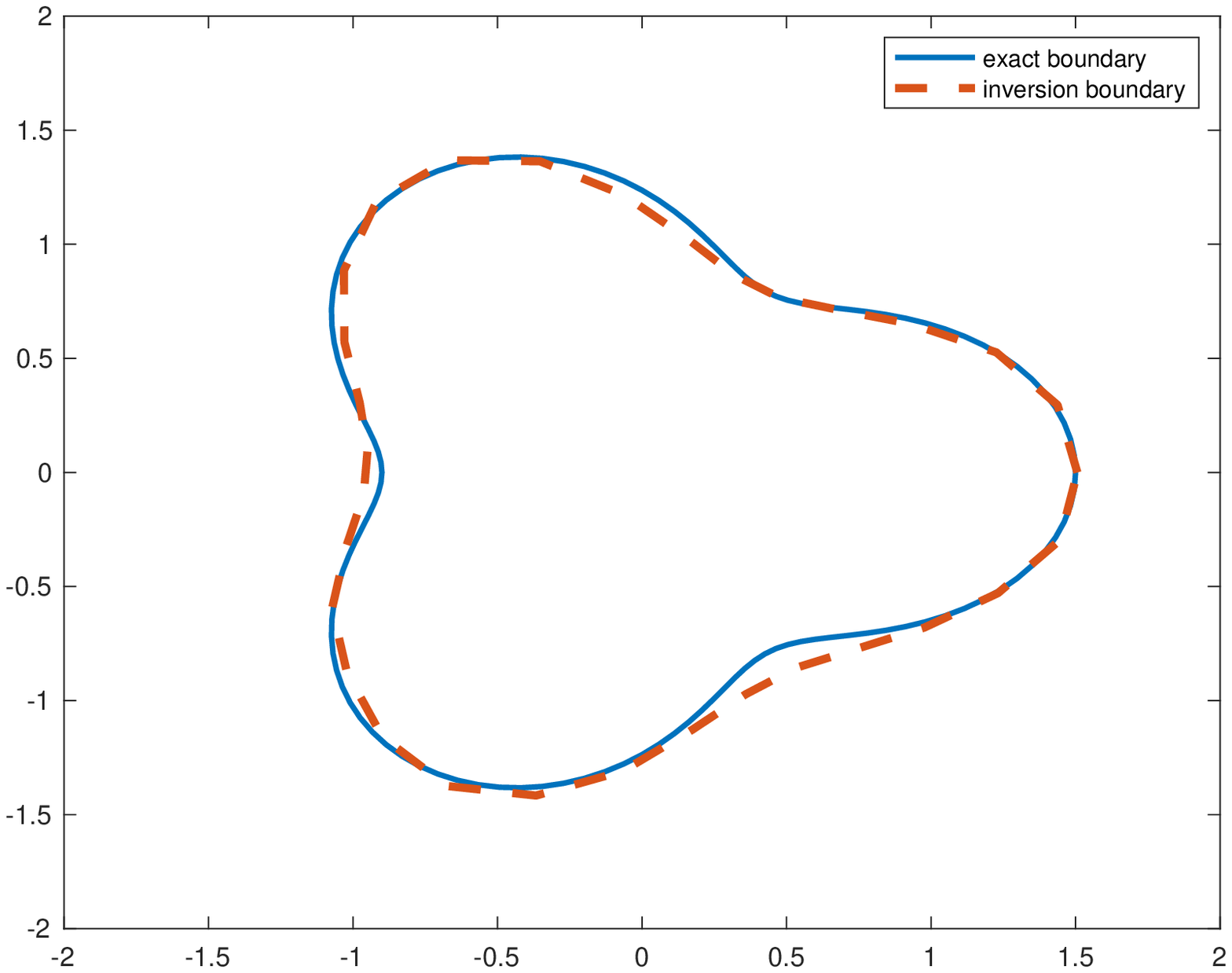}}\\
    \subfigure[Bean]{\includegraphics[width=4cm]{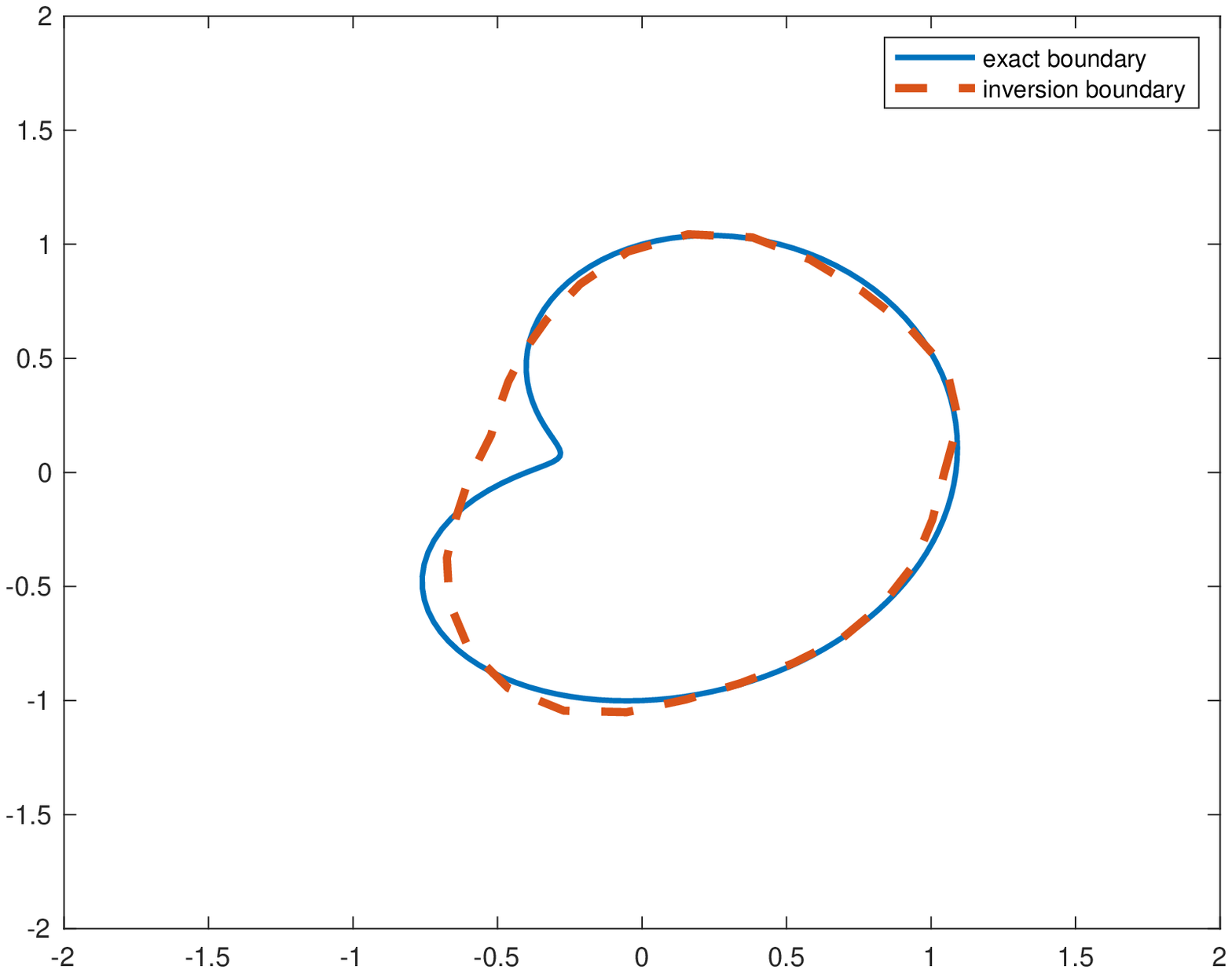}}
     \subfigure[Threelobes]{\includegraphics[width=4cm]{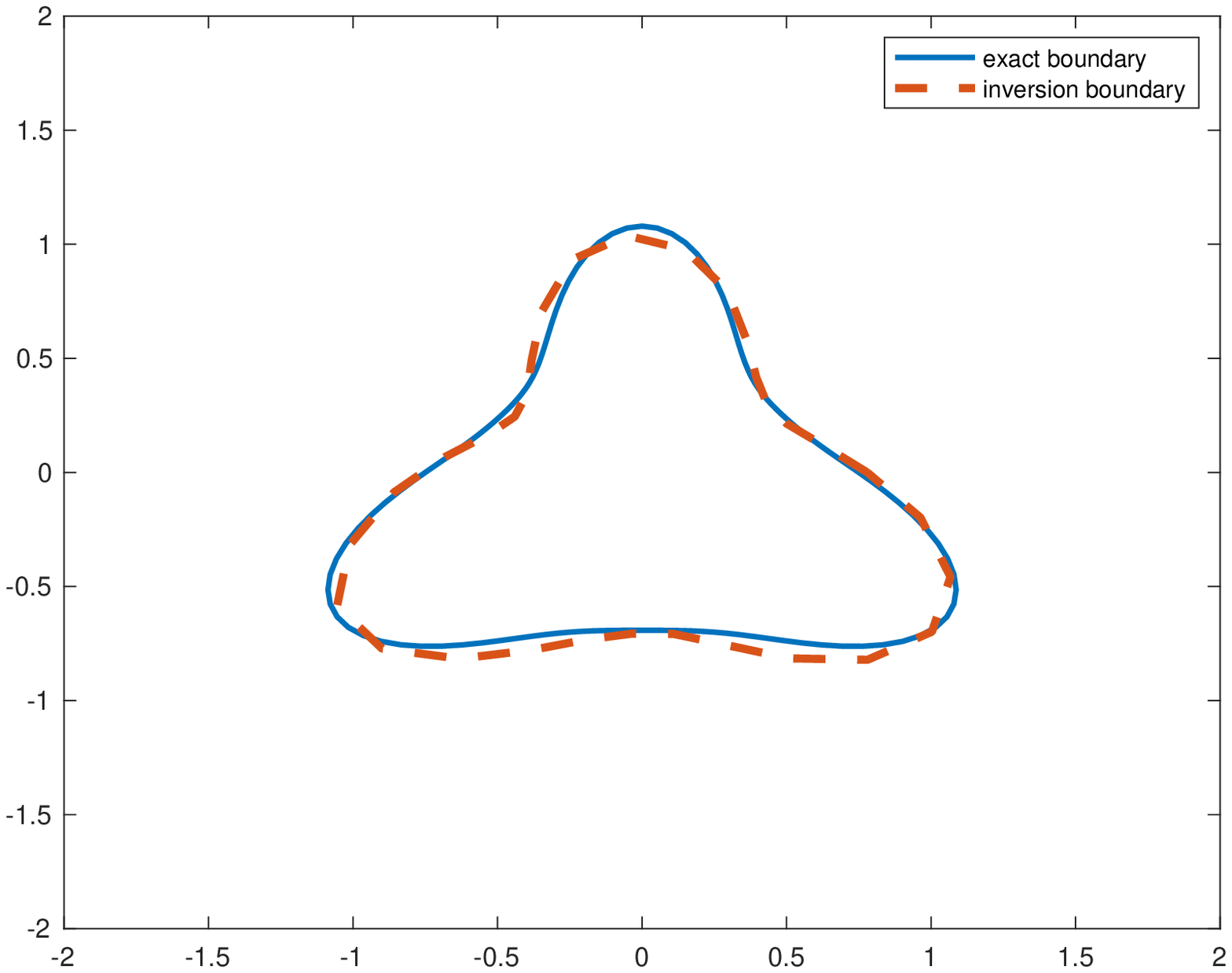}}\\
      \subfigure[Star]{\includegraphics[width=4cm]{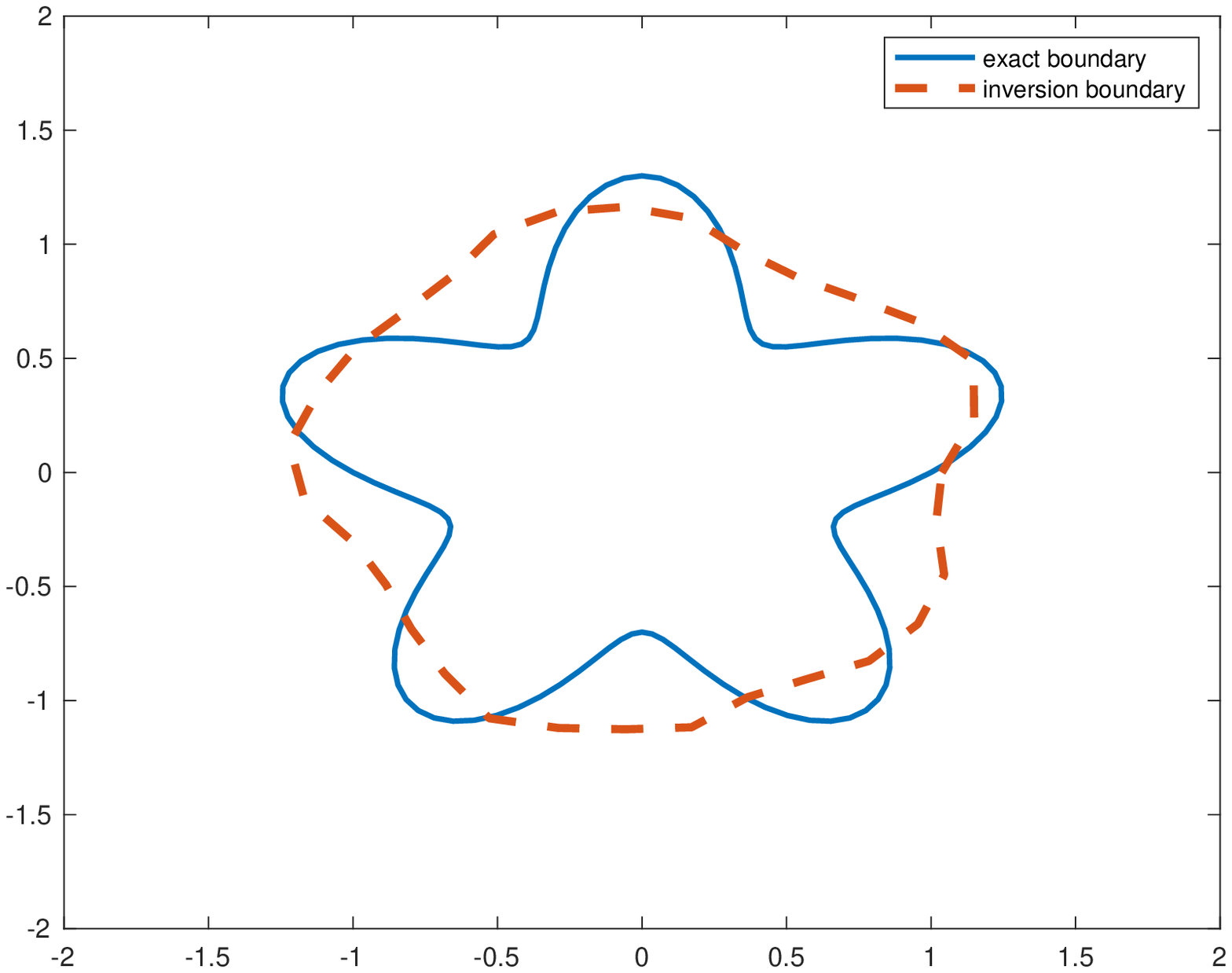}}
       \subfigure[Cloverleaf]{\includegraphics[width=4cm]{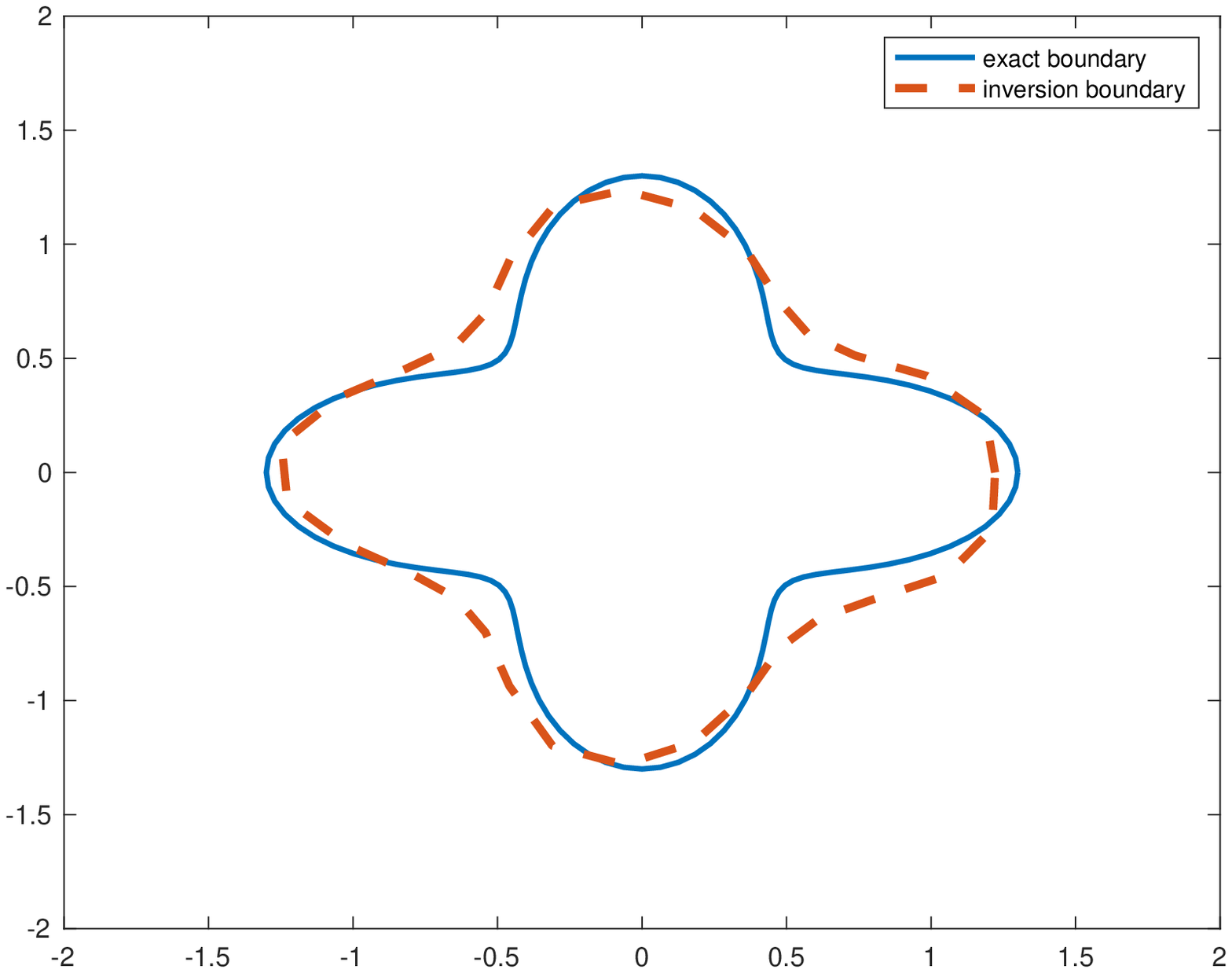}}\\
          \subfigure[Peanut]{\includegraphics[width=4cm]{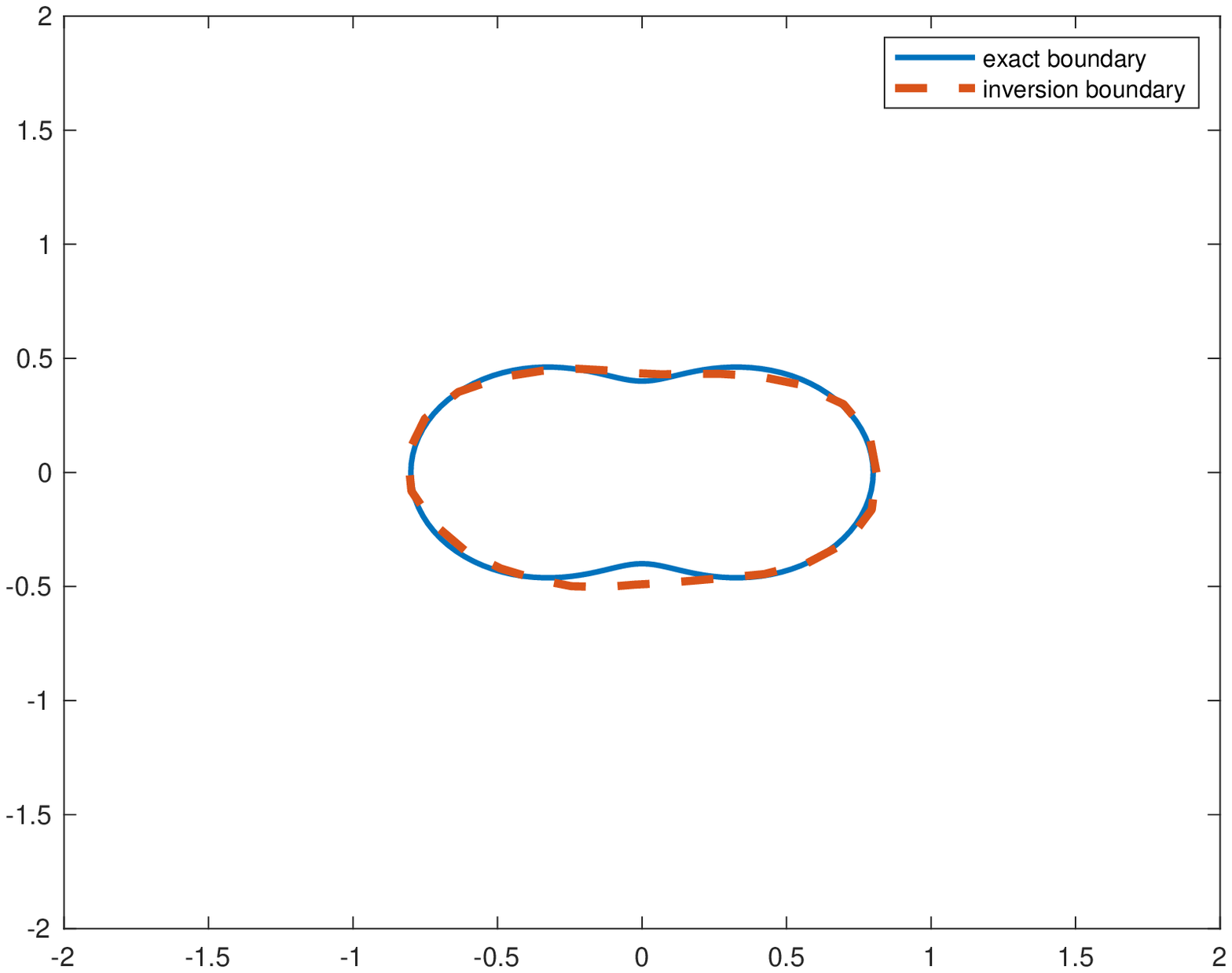}}
       \subfigure[Drop]{\includegraphics[width=4cm]{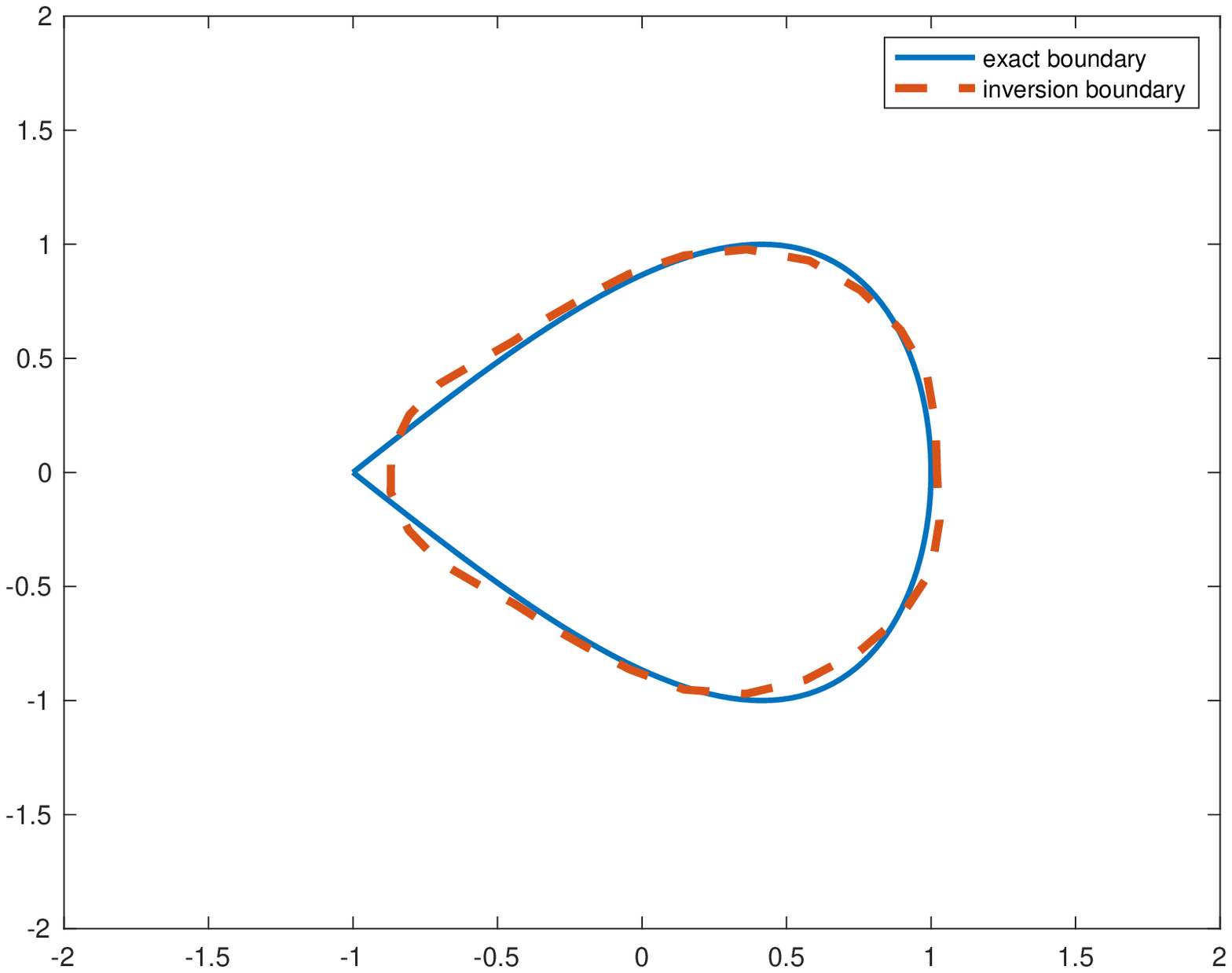}}
   \caption{Reconstructions for different shapes using $u^\infty(\hat{x}, d)$, $(\hat{x}, d)\in \gamma^{\rm o}_{2}\times\gamma^{\rm i}_{2}$ with $\eta_1=1\%, \,\eta_2=1\%$.}
    \label{resul3}
\end{figure}

\begin{figure}
\captionsetup[subfigure]{labelformat=empty}
    \centering
    \subfigure[Kite]{\includegraphics[width=4cm]{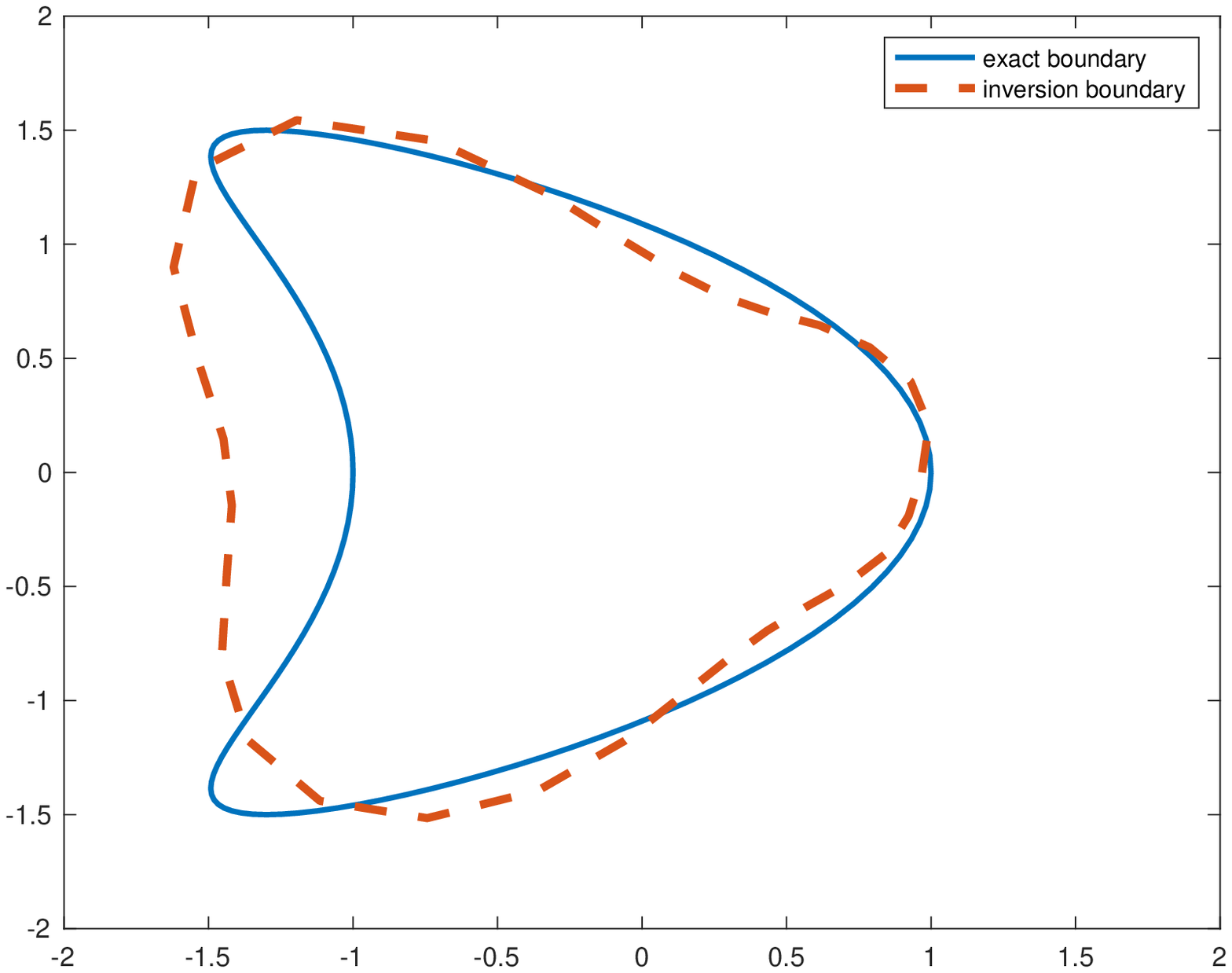}} 
    \subfigure[Roundrect]{\includegraphics[width=4cm]{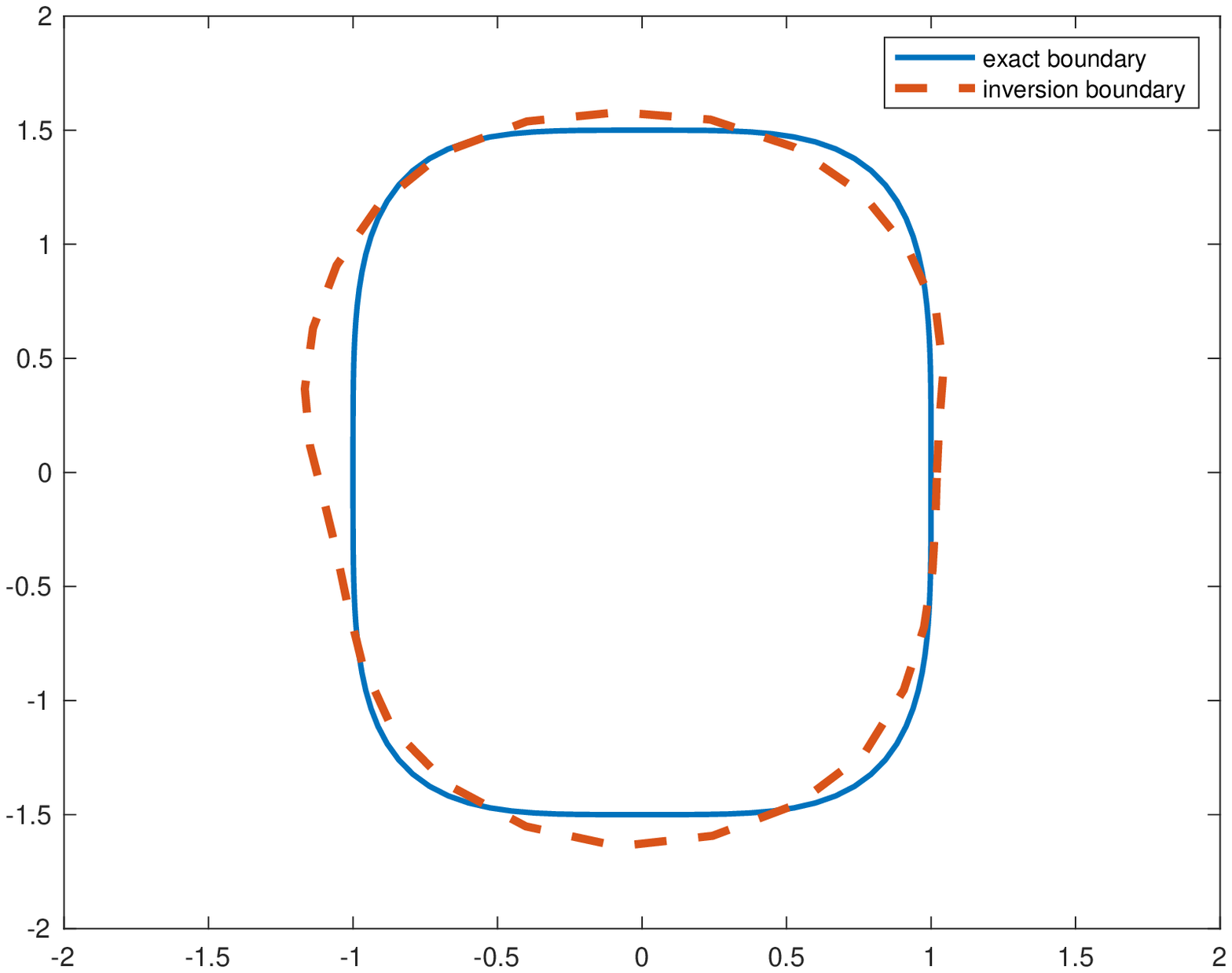}} \\
    \subfigure[Pear]{\includegraphics[width=4cm]{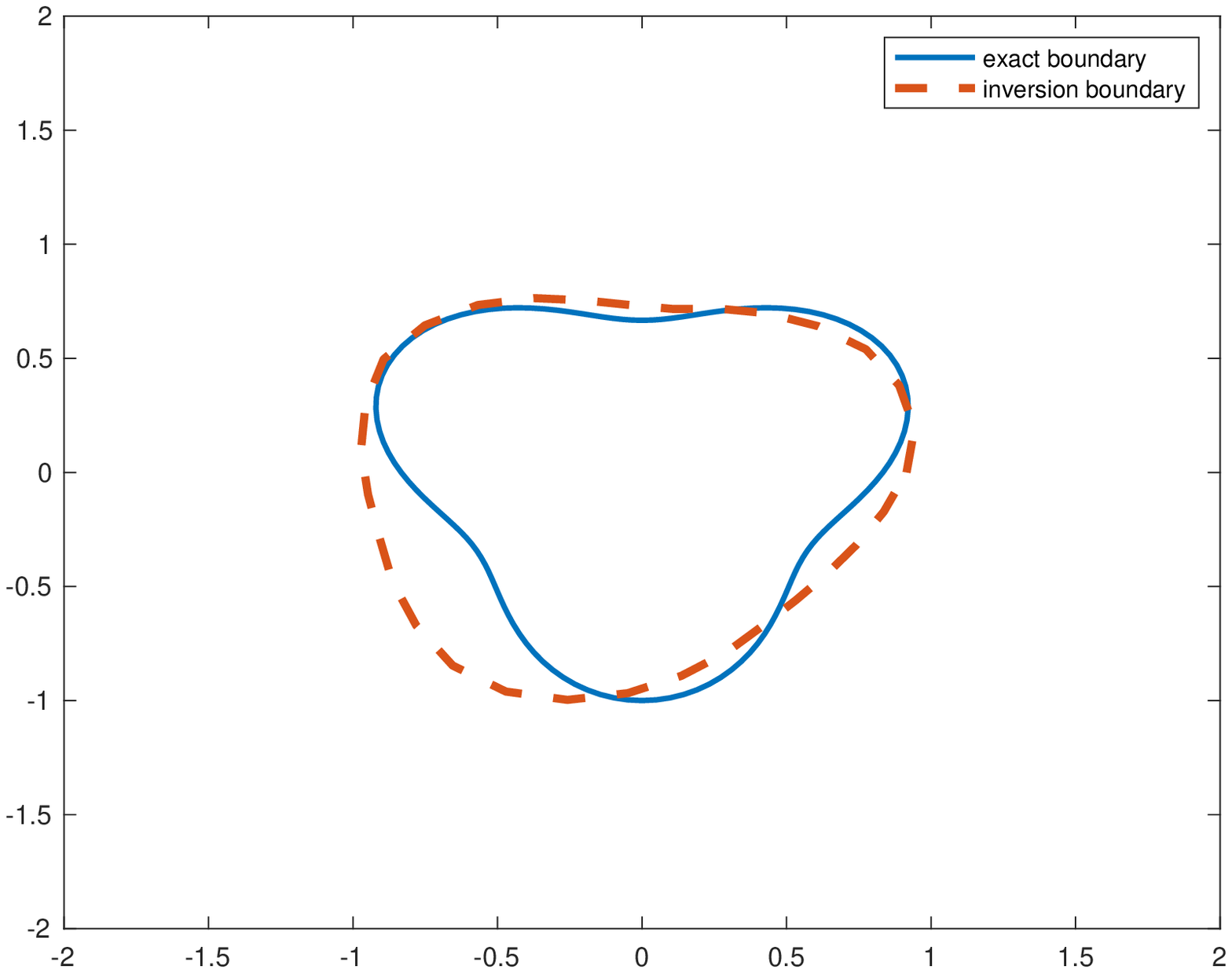}} 
 \subfigure[Acorn]{\includegraphics[width=4cm]{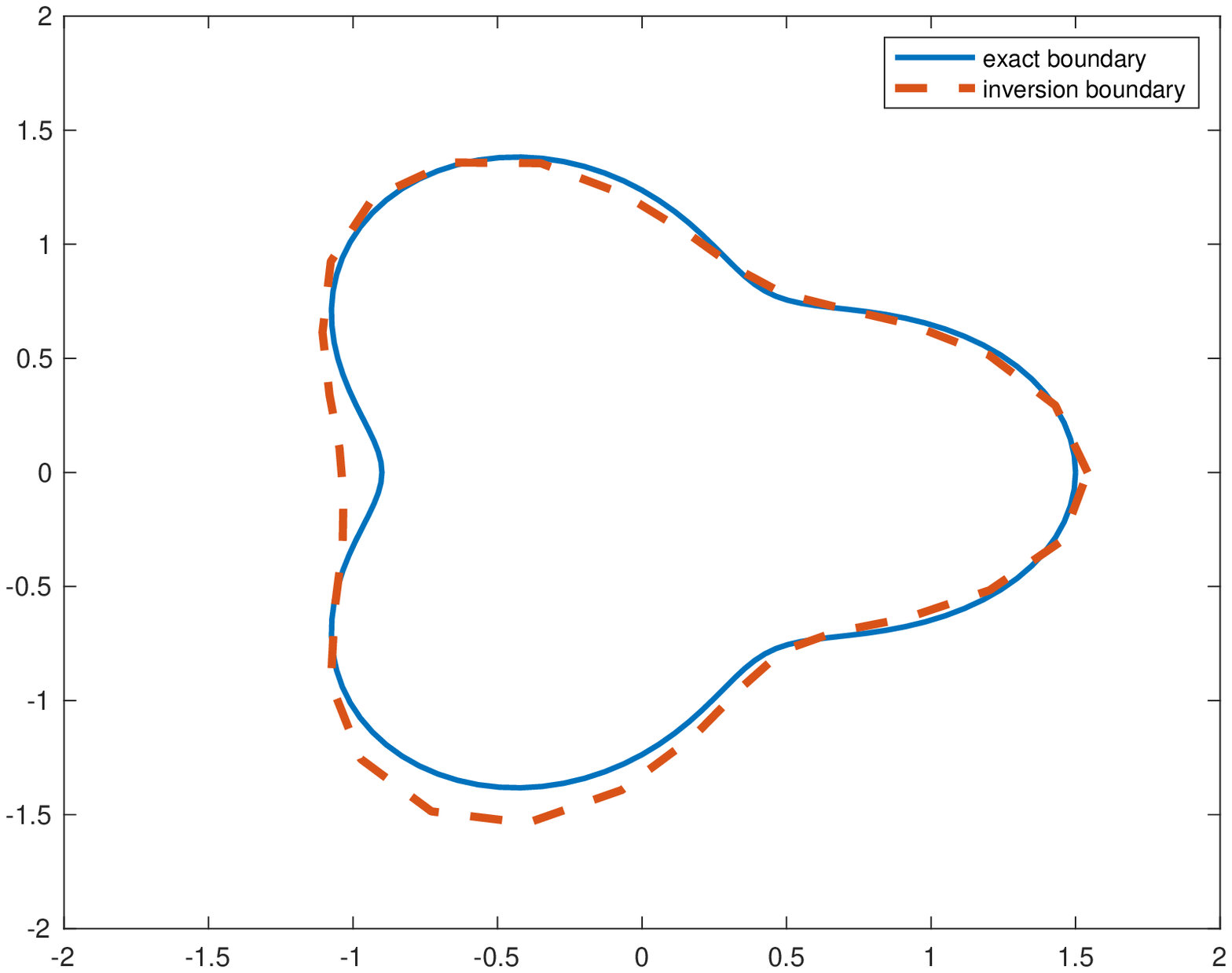}}\\
    \subfigure[Bean]{\includegraphics[width=4cm]{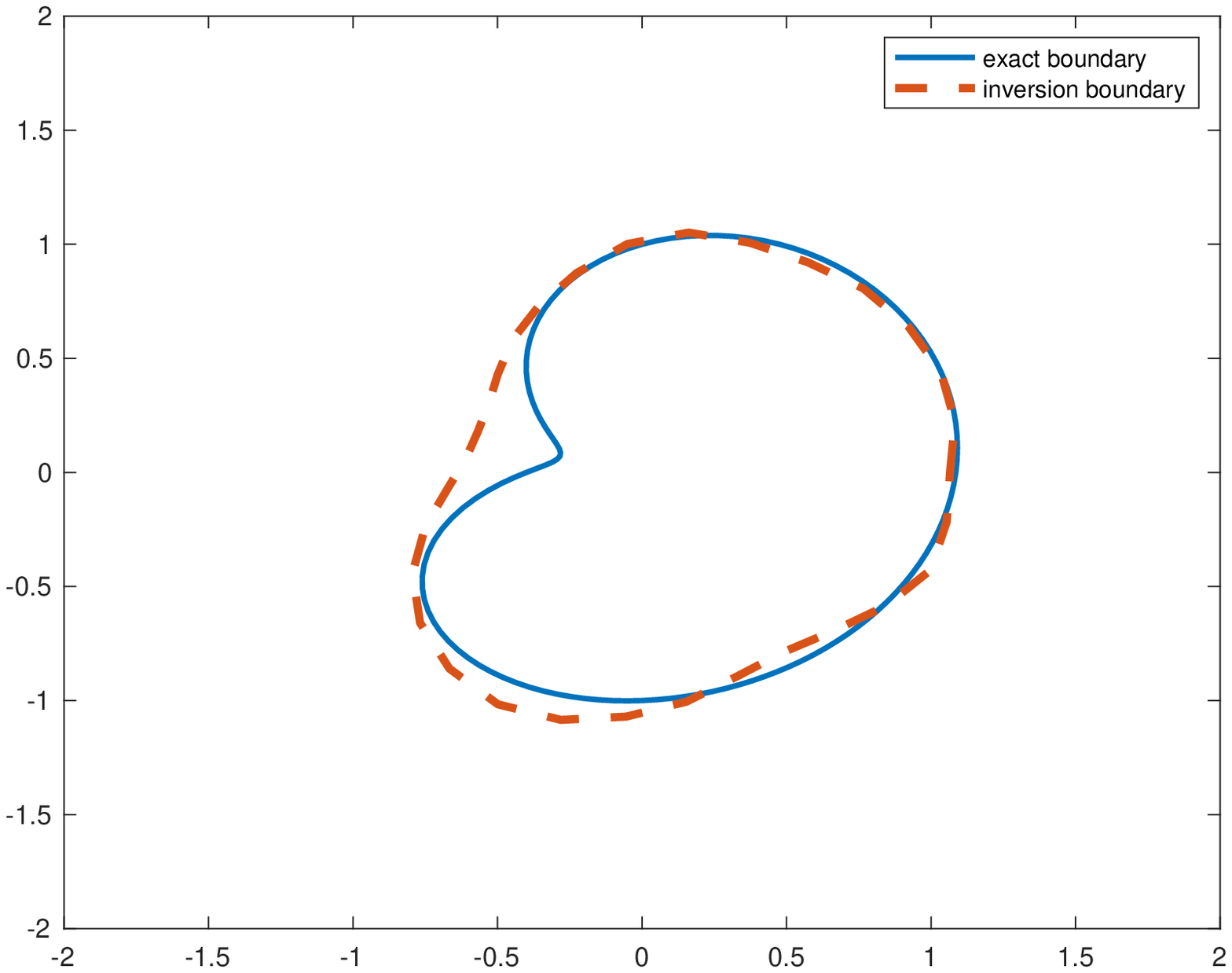}}
     \subfigure[Threelobes]{\includegraphics[width=4cm]{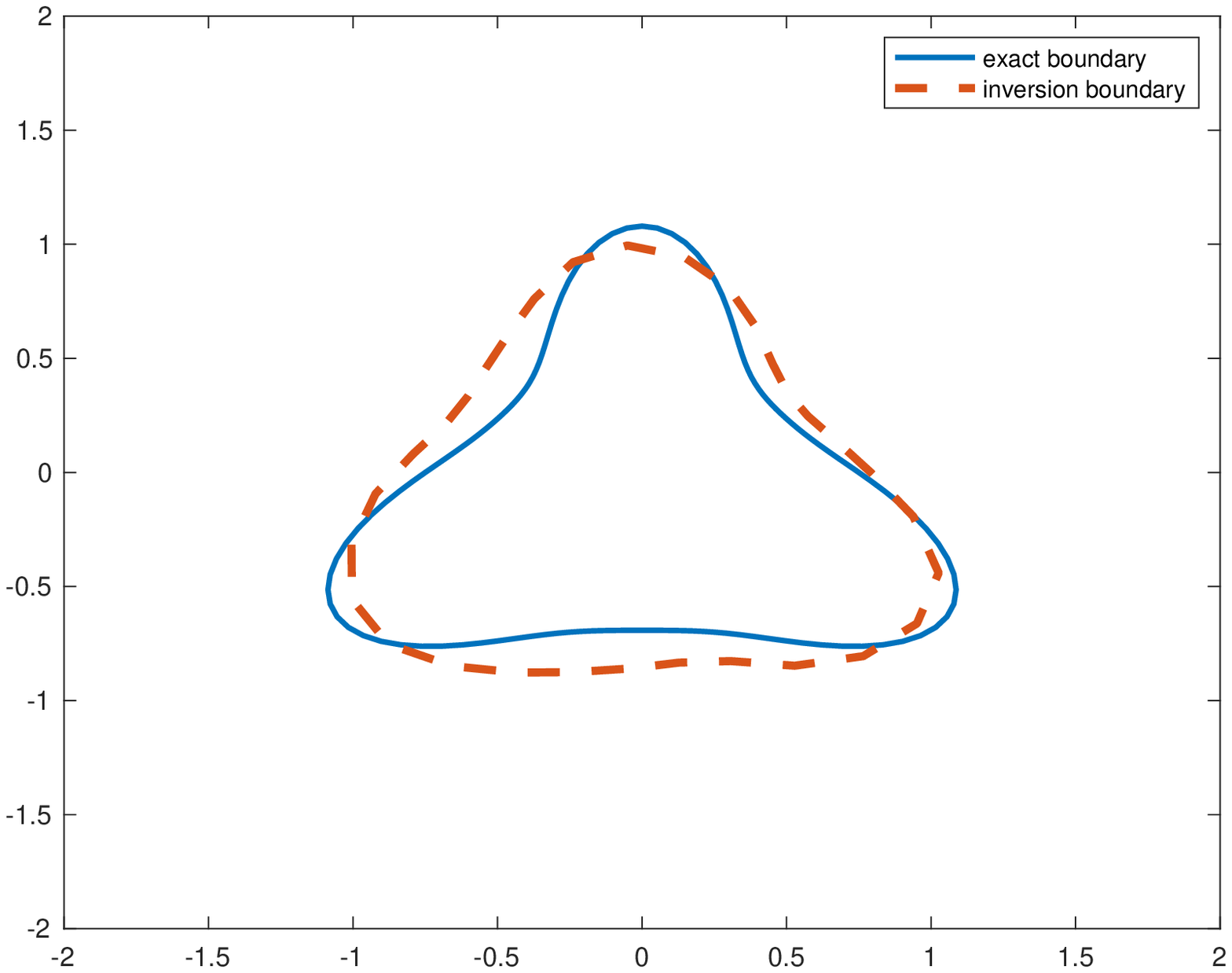}}\\
      \subfigure[Star]{\includegraphics[width=4cm]{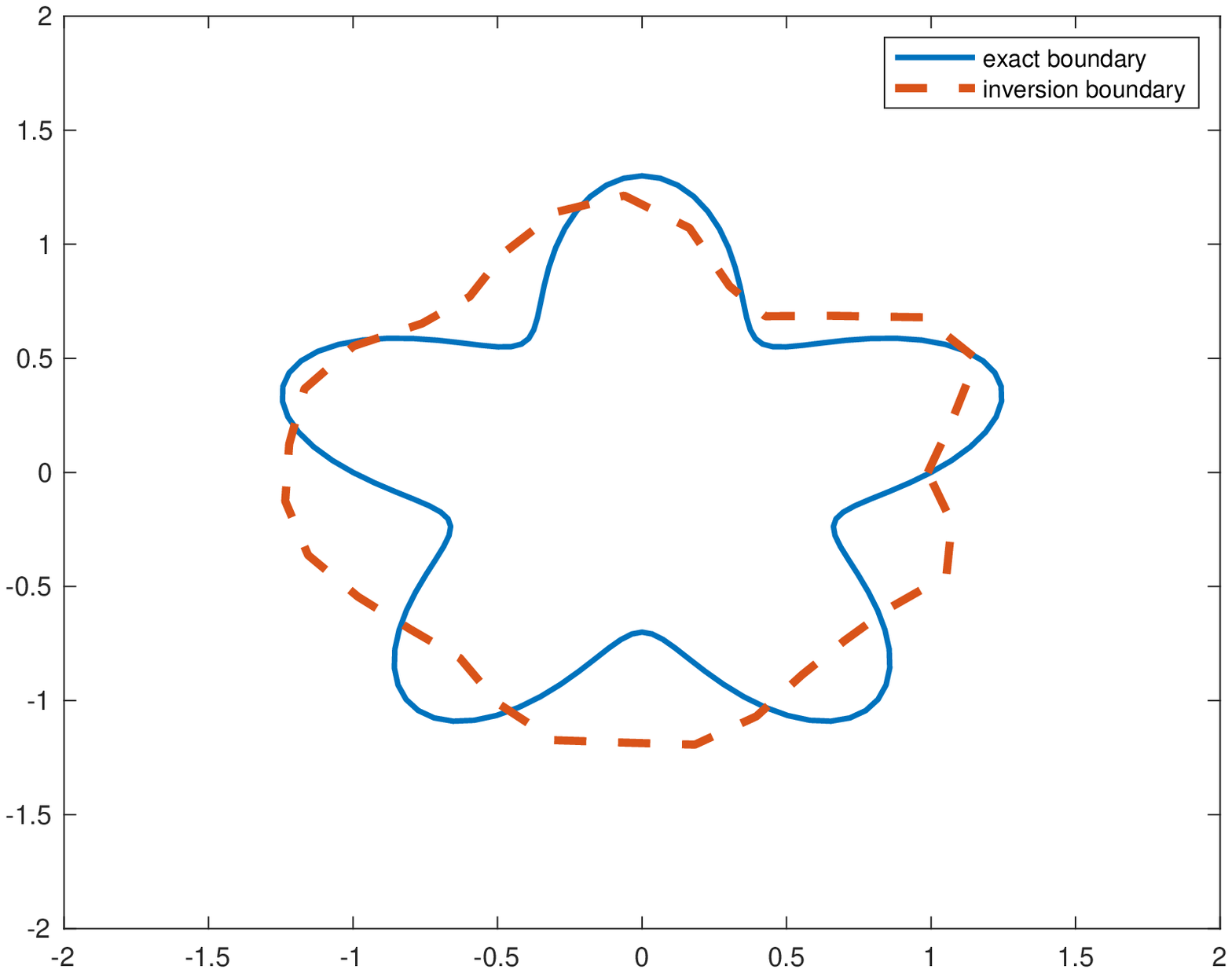}}
       \subfigure[Cloverleaf]{\includegraphics[width=4cm]{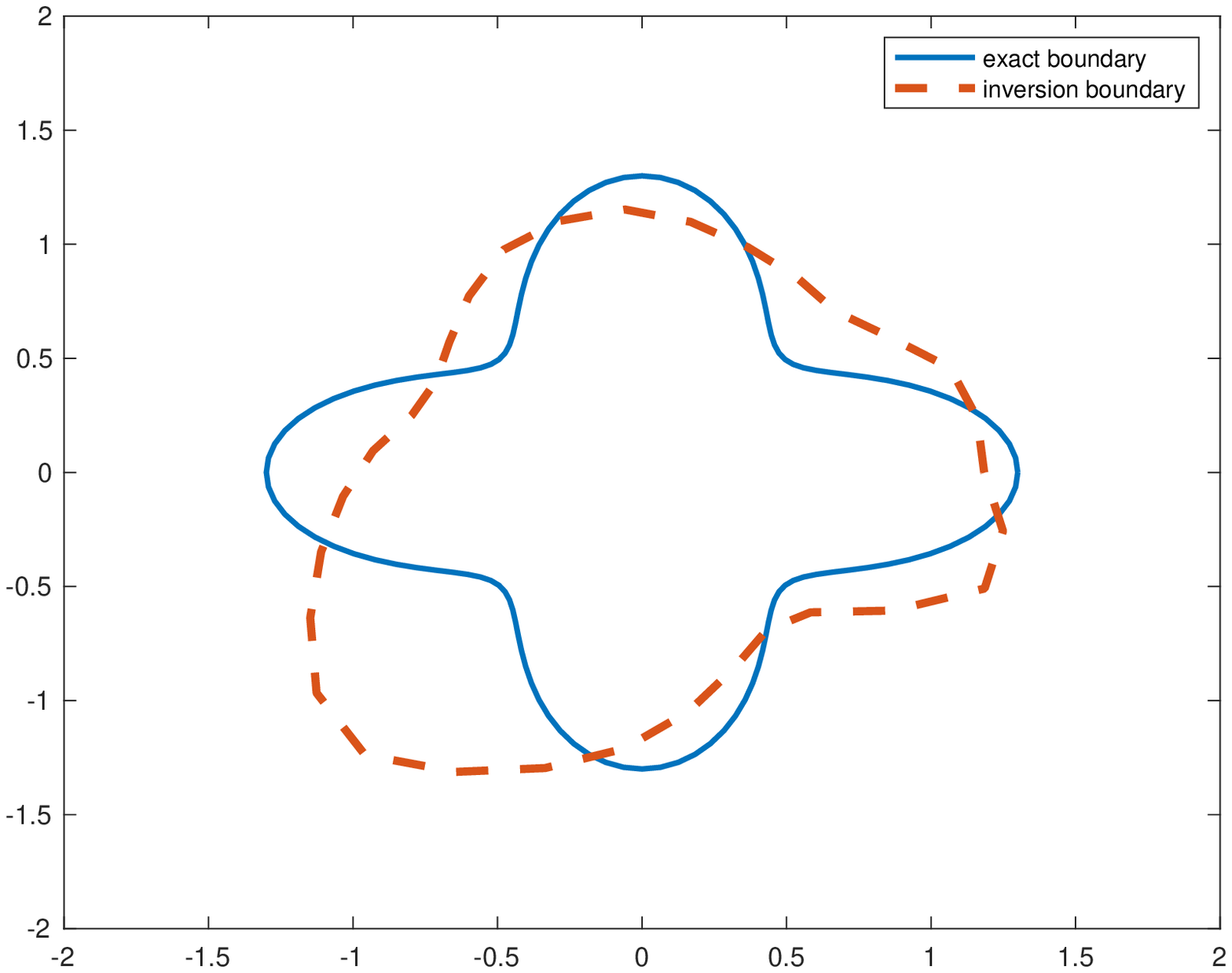}}\\
          \subfigure[Peanut]{\includegraphics[width=4cm]{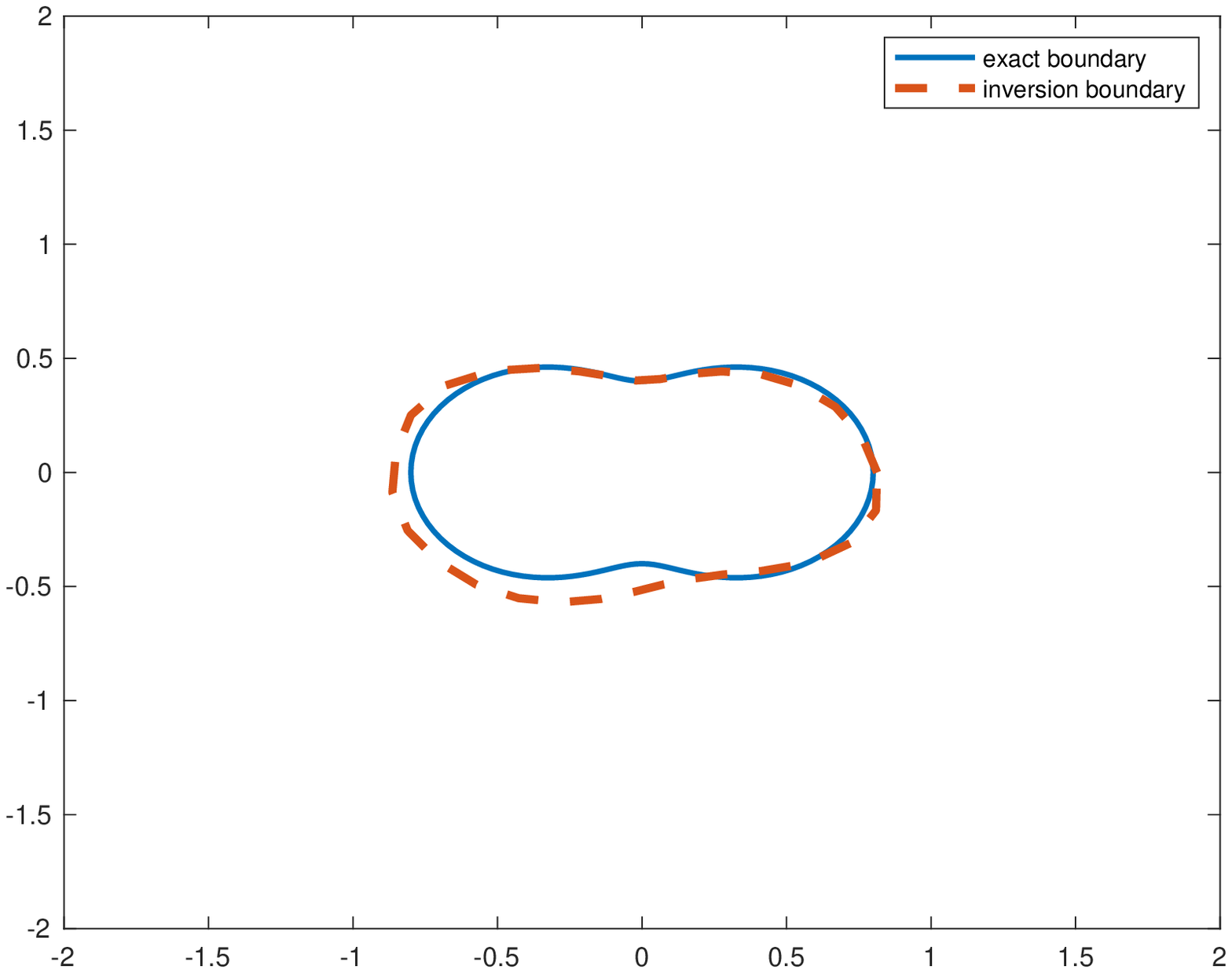}}
       \subfigure[Drop]{\includegraphics[width=4cm]{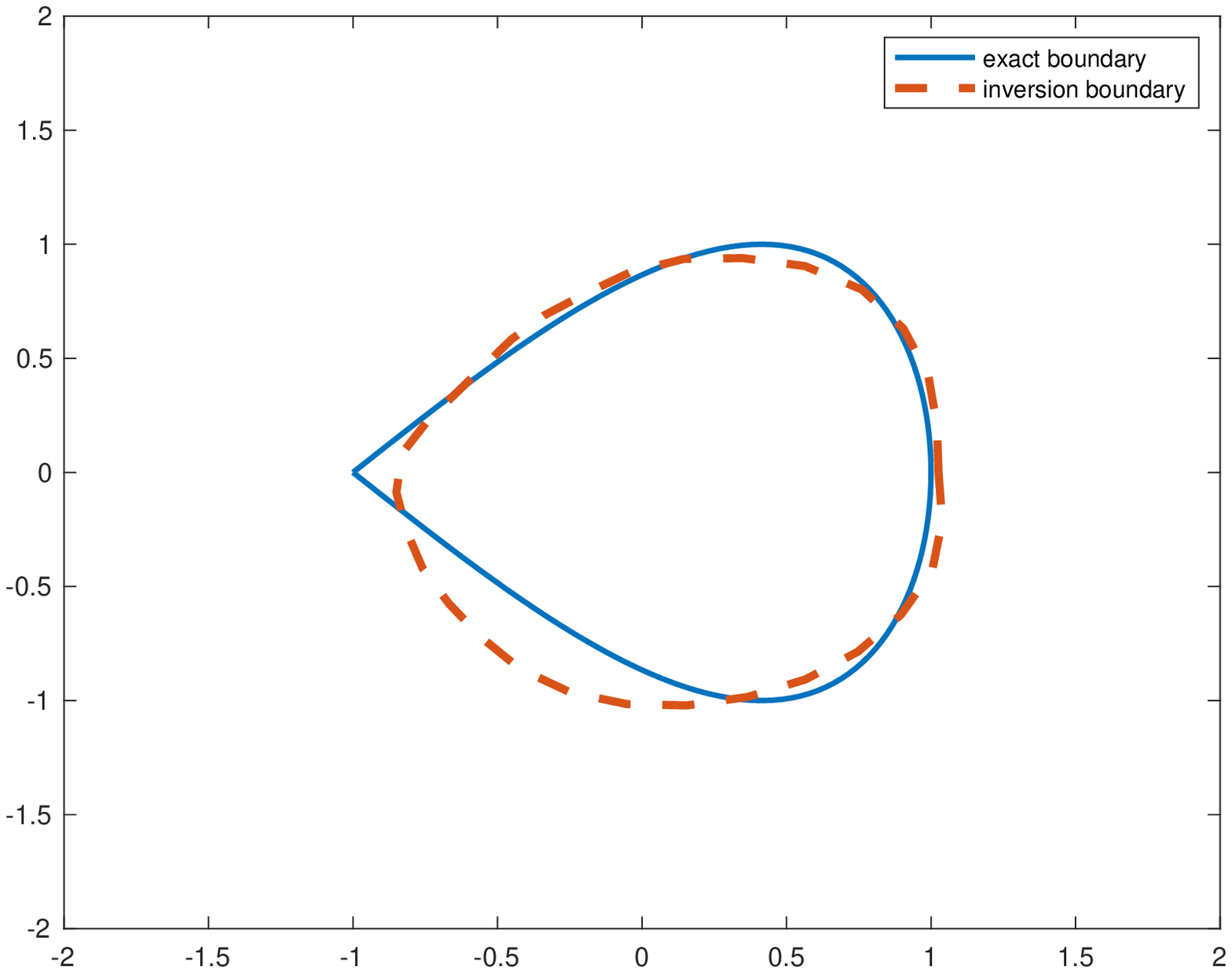}}
   \caption{Reconstructions for different shapes using $u^\infty(\hat{x}, d)$, $(\hat{x}, d)\in \gamma^{\rm o}_{3}\times\gamma^{\rm i}_{2}$ with $\eta_1=1\%, \,\eta_2=1\%$.}
    \label{resul4}
\end{figure}

The limited aperture cases are then considered with the same noise configuration in the full aperture cases. In Fig. \ref{resul3},  we show the numerical reconstructions using observation aperture $\gamma^{\rm o}=\gamma_2^{\rm o}$ with two incident waves, i.e.,  $\gamma^{\rm o}\times\gamma^{\rm i}=\gamma_2^{\rm o}\times\gamma_{2}^{\rm i}$. With smaller observation aperture $\gamma^{\rm o}=\gamma_3^{\rm o}$, Fig. \ref{resul4} shows the reconstructions for two incident waves $\gamma^{\rm i}=\gamma_{2}^{\rm i}$. The results show that the reconstructions are satisfactory for more observation data. When observation aperture becomes small, the reconstructions get worse.

\begin{figure}
\centering
\subfigure[]{
\begin{tikzpicture}
\coordinate (A) at (0,0);	
\coordinate (B) at (0,0.4);

\coordinate (C) at (0.4,0); 

\coordinate (D) at (0.4,0.4);
\coordinate (E) at (-0.4,0.4);
\coordinate (F) at (-0.4,-0.4);
\coordinate (G) at (-0.4,0);
\coordinate (H) at (0,-0.4);
\coordinate (I) at (0.4,-0.4);

\fill[yellow] (A) circle(0.1);
\foreach \n in {B, C, D, E, F, G, H, I}
	\fill[red] (\n) circle(0.1);
\draw[->, thick] (-2,0)--(2,0); 
\draw[->, thick] (0,-2.4)--(0,2);
\draw[ultra thick,variable=\t, domain=0:359.8,smooth] plot({cos(\t)+0.65*cos(2*\t)-0.65}, {1.5*sin(\t)});
\end{tikzpicture}
}
\subfigure[]{\includegraphics[width=6cm]{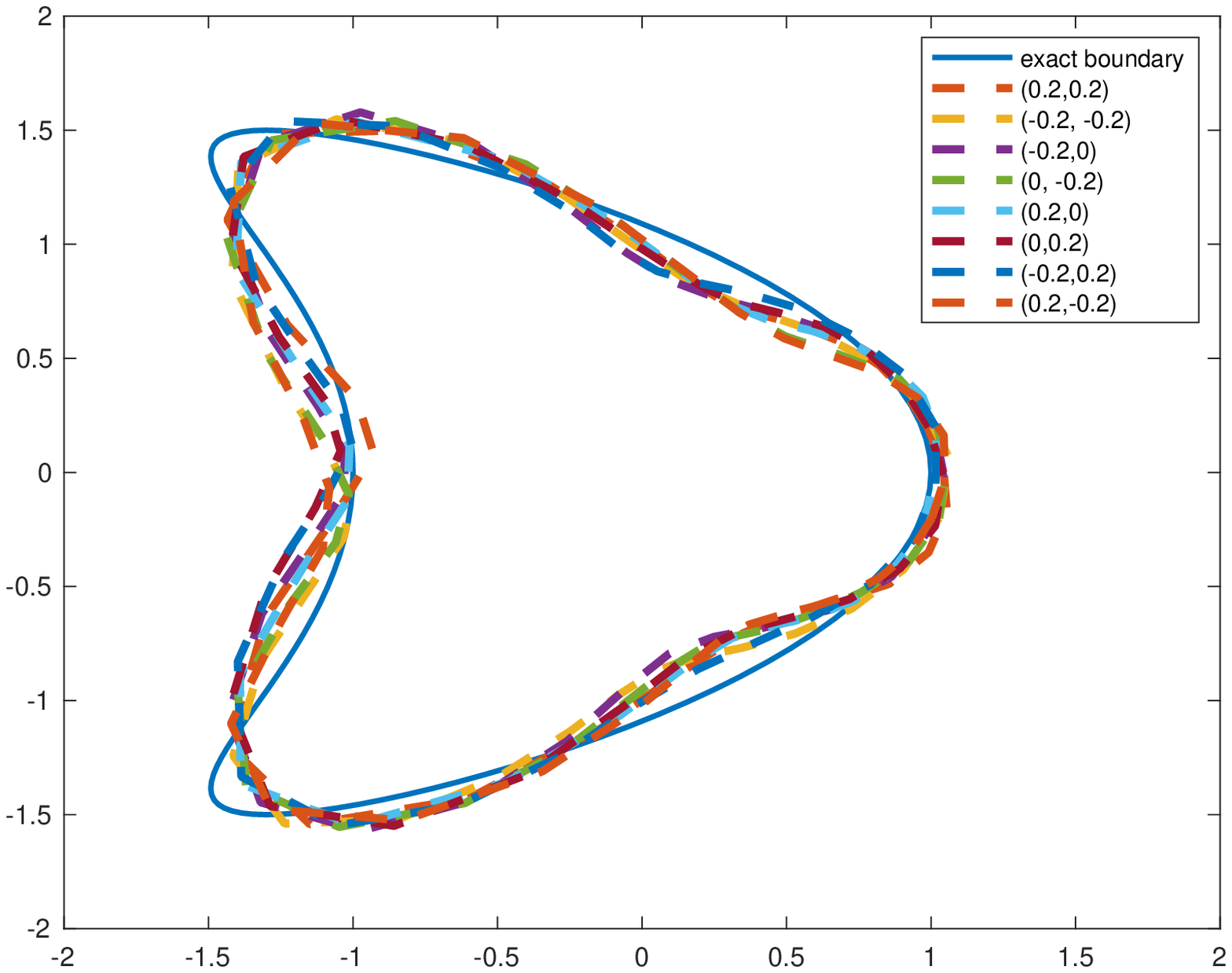}}
\caption{(a) Center points: Yellow dot is the true center position; Red dots are the center positions used in the reconstruction; (b) Reconstructions for Kite shape  using $8$ different centers with $(\hat{x}, d)\in \gamma^{\rm o}_{1}\times\gamma^{\rm i}_{2}$, $\eta_1=1\%, \,\eta_2=1\%$.}
\label{resul5}
\end{figure}

We also check the numerical accuracy using inaccurate center positions $(x_{\rm c}, y_{\rm c})$ in the reconstruction process. 
We use  the following $8$ different centers, (see Fig. \ref{resul5} (a))
$$(\pm 0.2, \pm 0.2), (0, \pm 0.2), (\pm 0.2, 0),$$
to reconstruct the Kite shape. The reconstructions are displayed In Fig. \ref{resul5} (b). From the displayed results, we can conclude that even the center is not exact, the reconstruction has not been seriously affected. This means that if we can use some method, e.g., the extended sampling method \cite{Li}, to find the approximated center first, the proposed method will provide satisfactory reconstruction. 


\section{Conclusions}
We discuss the application of the Bayesian method for the limited aperture inverse scattering problem. A novel total variation prior is proposed for the shape parameterization representation. Numerical examples show that the proposed method with the prior can yield satisfactory reconstructions even when the measurement data is quite limited. 

\bibliographystyle{plain}\bibliography{ref}

\end{document}